\documentclass[preprint,12pt]{elsarticle}
\usepackage[utf8]{inputenc}
\usepackage{textcomp}
\usepackage{pstricks}
\usepackage{amssymb}
\usepackage{pifont}
\usepackage{graphicx}
\usepackage[version=4]{mhchem}
\usepackage{siunitx}
\usepackage{longtable,tabularx}
\usepackage[numbers]{natbib}
\usepackage{hyperref}
\usepackage[capitalise]{cleveref}
\usepackage{amsfonts,amsmath,mathtools,amssymb}
\usepackage{setspace}
\usepackage{placeins}
\usepackage[colorinlistoftodos]{todonotes}
\usepackage{booktabs}
\usepackage{float}
\usepackage{algpseudocode}
\usepackage{algorithm}
\usepackage{subfigure}
\usepackage{lineno}
\usepackage{tikz}
\usepackage{multirow}
\usepackage{rotating}
\usepackage{appendix}
\usepackage{makecell}
\usetikzlibrary{arrows,calc,through,backgrounds,matrix,decorations.pathmorphing}
\usepackage{soul}
\usepackage{rotating}
\sethlcolor{yellow}
\soulregister\ref7
\soulregister\cref7
\soulregister\cite7
\soulregister\citet7
\soulregister\eqref7
\usepackage{multirow}

\newcommand{\overbar}[1]{\mkern1.5mu\overline{\mkern-1.5mu#1\mkern-1.5mu\mkern1.5mu}}
\newcommand{\f}[2]{\frac{#1}{#2}}
\newcommand{\mb}[1]{\mathbf{#1}}
\newcommand{\tr}[1]{\mathrm{Tr}\left({#1}\right)}

\DeclareMathAlphabet\mathbfcal{OMS}{cmsy}{b}{n}
\renewcommand{\d}{\mathop{}\!\mathrm{d}} % total derivative
\newcommand{\p}{\partial}
\crefname{equation}{eq.}{eqs.}  % For lowercase "eq." and "eqs."
\Crefname{equation}{Eq.}{Eqs.}  % For capitalized "Eq." and "Eqs."

\Crefname{appendix}{}{}
\newcommand{\vop}{\mathrm{vec}}

\definecolor{rev1}{HTML}{FFFFFF}%{FF999A} % red
\definecolor{rev2}{HTML}{FFFFFF}%{F3F298} % yellow
\definecolor{rev3}{HTML}{FFFFFF}%{B2E0AE} % green
\definecolor{other}{HTML}{FFFFFF}%{C8C7FF} % blue
\DeclareRobustCommand{\hlone}[1]{{\sethlcolor{rev1}\hl{#1}}}
\DeclareRobustCommand{\hltwo}[1]{{\sethlcolor{rev2}\hl{#1}}}
\DeclareRobustCommand{\hlthree}[1]{{\sethlcolor{rev3}\hl{#1}}}
\DeclareRobustCommand{\hlo}[1]{{\sethlcolor{other}\hl{#1}}}

\begin{document}

\title{Differentiable Singular Value Decomposition (SVD)}

\author[inst1]{Rohit Kanchi}
\author[inst1]{Sicheng He}

\affiliation[inst1]{organization={Department of Mechanical, Aerospace, and Biomedical Engineering, University of Tennessee},
    city={Knoxville},
    postcode={TN 37996},
    country={USA}}

\begin{frontmatter}
\begin{abstract}
Singular value decomposition (SVD) is widely used in modal analysis, such as proper orthogonal decomposition (POD) and resolvent analysis, to extract key features from imaginary problems.
SVD derivatives need to be computed efficiently to enable the large-scale design optimization.
However, for a general imaginary matrix, no method can accurately compute this derivative to machine precision and remain scalable with respect to the number of design variables without requiring the all of the singular variables.
We propose two algorithms to efficiently compute this derivative based on the adjoint method and reverse automatic differentiation (RAD) and RAD-based singular value derivative formula.
Differentiation results for each method proposed were compared with FD results for one square and one tall rectangular matrix example and matched with the FD results to about 5--7 digits.
Finally, we demonstrate the scalability of the proposed method by calculating the derivatives of singular values with respect to the snapshot matrix derived from the POD of a large dataset for a laminar-turbulent transitional flow over a flat plate, sourced from the John Hopkins turbulence database (JHTDB). 
\end{abstract}
\end{frontmatter}

\section{Introduction}\label{sec:intro}
% Topic: SVD application
Singular variable decomposition (SVD) is a fundamental linear algebra technique widely used in engineering applications~\cite{Gilbert2022,eckart1939}.
SVD is applied in varied fields such as numerical optimization~\cite{candes2012a}, mechanical design~\cite{gwak2004a,Khalil2007,Sarkar2009a,lucas2012a,Zhao2016b}, fluid dynamics~\cite{Li2023a,diaz2024a}, noise-reduction~\cite{shin2003a,Pan2019,Li2023}, image processing~\cite{sadek2012a}, and finance~\cite{wang2017b}.
Applications of SVD in engineering optimization and analysis problems include aerostructural optimization~\cite{lucas2011a,choi2020a,poole2019a}, aerodynamic shape optimization~\cite{Li2018a,li2019g,yao2020a,Li2021b,Wu2023a,chen2024a},\hlone{ vibration analysis~\cite{Liu2024,Li2022d}} and pure structural optimization problems~\cite{liu1997a,lee1999a,ersoy2002a,ersoy2013a}.
Many modal analysis methods are based on SVD, e.g., the resolvent analysis~\cite{Trefethen1993}, and dynamic mode decomposition (DMD)~\cite{SCHMID2010}.
For more details about SVD in modal analysis, see the review papers by Taira et al.~\cite{Taira2019a,Taira2020}.

%Topic: SVD computation approaches (for gov equn)
There are several routes to compute SVD of a given matrix (refer to~\citet[Chapter~8.6]{Golub1996} for more details).
A popular approach is to reduce the problem to an eigenvalue problem (EVP)~\cite{Golub1996,Trefethen1997a}
Once that is done, numerical methods such as Golub-Kahan bi-diagonalization~\cite{Golub1996} and the Jacobi SVD algorithm~\cite{Demmel1992} can be used to compute the SVD. 
Iterative but approximating methods such as Lanczos or Krylov subspace methods~\cite{Saad2011}, randomized SVD~\cite{Halko2011}, and truncated SVD have also been developed to reduce computational costs, particularly for large-scale data or high-dimensional systems. 

These advancements are crucial when differentiation of the SVD is needed, as the computational cost and numerical stability of differentiation depend heavily on the underlying method used to compute the SVD.
This motivates the need to carefully evaluate SVD computation strategies for their efficiency and compatibility with different differentiation frameworks.
For a more thorough review of SVD history and algorithms, especially the relationship between EVP and SVD, we direct readers to the review by Zhang~\cite{Zhang2015d} and Lecture 31 of Trefethen and Bau~\cite{Trefethen1997a}.

% Topic: GMM and SEM (description only, no derivative method mention yet)
Two approaches in the relationship between the EVP of a symmetric matrix and the SVD were leveraged in the current study.
In the first approach, we perform eigendecomposition of the product of the matrix and its Hermitian transpose (if real, just transpose), and the reverse product~\cite{Golub1996,Trefethen1997a,Salgado2022}
The resulting matrices from this pair of products are the Gram matrices~\cite{Schwerdtfeger1961,Trefethen1997a,Makkonen2023} of the matrix under consideration whose SVD is sought, or in other words, the target matrix.
The eigendecomposition of the Gram matrix resulting from the product of the target matrix and its Hermitian transpose gives the singular values and the left singular vectors of the SVD.
Here, the eigenvectors are the left singular vectors, and the square roots of the eigenvalues are the singular values.
Similarly, the eigendecomposition of the Gram matrix from the product of the reverse product of the target matrix and its Hermitian transpose yields the singular values and the right singular vectors of the target matrix~\cite{Golub1996}.
We developed governing equations for this approach based on the approach taken by He et al.~\cite{He2023}.
This approach is useful when the matrix whose SVD is sought is sparse.
If it is dense, then forming the left or right Gram matrices is costly, which motivated the second approach.

In the second approach, we leverage the relationship between the SVD of the target matrix and the EVP of its symmetric embedding~\cite{Golub1996,Ragnarsson2013}.
The relationship is such that the singular value and its negative counterpart are eigenvalues of the symmetric embedding.
The eigenvectors of the symmetric embedding contain the left and right singular vectors of the target matrix.
We developed a set of SVD governing equations in the current manuscript based on this symmetric embedding which can be solved using an alternating least squares approach with Rayleigh quotient through successive deflation~\cite{Golub1996,Trefethen1997a,Hansen1998}.%RK-HS-: I have a code that does this and I have pushed it to repo. You may check it.

%Topic: Repeated sing vals
The derivative computation methods proposed in the current manuscript are not limited to only the EVP-SVD approach to compute the singular variables.
So long as the singular variables computed through any algorithm satisfy the governing equations for either aforementioned case, the proposed derivative algorithms will work.
However, there arise problems whose solutions might involve nearly similar or duplicated singular values~\cite{Golub1996,Stewart1998,Angelova2024}.
Since the singular values are the square roots of the eigenvalues of the Gram matrices, the Gram matrices EVPs yield repeated eigenvalues in the case of duplicate singular values.
Typically, such duplicated or close eigenvalues occur due to some spatial symmetry in the problem whose eigen-solution is sought~\cite{Davies2004,Lin2020}.
It is a fact that the changes of the order of $\epsilon$ in a matrix can alter its singular subspace by $\delta / \epsilon$ where $\delta$ is a measure of how far apart the singular values are from each other~(refer to~\citet[Chapter~8.6.1]{Golub1996}).
Thus, if two singular values are close or the same, the ratio would tend to infinitely large values.
This means that the singular vector space is mathematically degenerate~\cite{Lin2020}.
Angelova and Petkov~\cite{Angelova2024} showed in a component-wise perturbation analysis that the SVD perturbation problem is only well-posed when the singular values are distinct.
The distinct singular values, therefore, ensure uniqueness in solution to the governing equations developed.
This case of duplicate singular values is, however, beyond the scope of the current manuscript.

% Topic: SVD derivative application
Once the required singular variables are obtained from either of the two EVP-SVD approaches mentioned, the derivative of the singular variables is sought, which has applications in different design optimization problems. 
For instance, Ersoy and Mugan~\cite{ersoy2002a} computed the derivatives of the singular values and vectors for the structural design derivative analysis.
Santiago et al.~\cite{santiago2021a} took the derivatives of the singular values of the shapes in the shape derivative as part of the shape optimization for photonic nanostructures.
Recently, Skene and Schmid~\cite{Skene2018a} performed a mixed forward and backward derivative analysis for the swirling M-flames, which involved computation of the resolvent with respect to the design parameters.
The resolvent analysis involved the singular value decomposition of the resolvent matrix, which gave the dominant singular value whose derivative was sought.

% Topic: General derivative computation method
Various methods are used for computing derivatives, including finite differences (FD), imaginary step (CS), algorithmic differentiation (AD), the direct method, and the adjoint method (refer to~\citet[Chapter~6]{Martins2022}). 
These methods differ in terms of their accuracy and computational efficiency. 
Efficiency-wise, each method is typically optimized for either scaling with the number of outputs (e.g., functions of interest such as eigenvalues and eigenvectors) or the number of inputs (e.g., design variables), but rarely both~\cite{Martins2013a},\cite[Chapter~6]{Martins2022}.

% Topic: FD/CS
FD is susceptible to truncation and subtraction errors, whereas CS avoids these issues with a sufficiently small step size and can achieve machine precision~\cite{Martins2003a}.
FD and CS are relatively easy to implement due to their black-box-like nature.  
However, these methods may not be practical for high-fidelity applications involving a large number of design variables. (For a detailed comparison of FD and the adjoint method, refer to \cite{Lyu2014f} and \cite[Chapter~6]{Martins2022} (Fig. 6.43)).
Further, their computational cost scales unfavourably with the number of inputs, with CS being more expensive because of the use of imaginary arithmetic.
For instance, with respect to differentiating SVD, Ersoy, and Mugan~\cite{ersoy2002a} computed the SVD derivatives using the FD approach, which lacks both accuracy and scalability with respect to design variables.
This motivates us to choose a method that does not scale with a number of inputs.

% Topic: AD
AD differentiates a program based on the systematic application of the chain rule.
Forward-AD (FAD) does this by applying the chain rule from inputs to the outputs, and reverse-AD (RAD) does so in the reverse order. 
Thus, the cost is proportional to the number of inputs when FAD is implemented and the number of outputs when RAD is used.
Currently, for derivative algorithms that do not scale with inputs, there are only RAD-based methods in the literature to compute the SVD derivative~\cite{Giles2000,townsend2016a,seeger2017a,wan2019a}.
The RAD formula proposed by Giles~\cite{Giles2000} works for only real values. 
The RAD form proposed by Townsend~\cite{townsend2016a} works for real values and works with reduced SVD.
Santiago et al.~\cite{santiago2021a} took the derivatives of the singular values for the shape optimization for photonic nanostructures by using the RAD formula given for real-valued SVD by Townsend~\cite{townsend2016a}.
Wan and Zhang~\cite{wan2019a} proposed a RAD formula for the full SVD that works with complex-valued inputs.
Finally, Seeger et al. ~\cite{seeger2017a} proposed a formula equivalent to Townsend's~\cite{townsend2016a}.
However, these methods require all of the singular variables during the derivative computation.
This increases the computational storage required, which motivates the proposal of a method that does not involve several large matrix-matrix products and computes the derivative more memory-efficiently.
Projection-based method was developed by He et al.~\cite{He2022a} to partially address this challenge in eigenvector derivative computation.

% Topic: direct/adjoint
Beyond explicit analytic methods like AD, implicit analytic methods like direct and adjoint approaches can also be used~\cite[Chapter~6]{Martins2022}. 
The efficiency of these methods depends on the relative number of inputs and outputs. 
The direct method is preferable when the number of inputs is smaller than the number of outputs. 
In contrast, the adjoint method is more efficient when the number of inputs exceeds the number of outputs.
Skene and Schmid~\cite{Skene2018a}, for instance, used the direct approach to compute the resolvent derivative.
However, this scales inefficiently with the number of inputs.
Further, there are no adjoint-based methods to compute the SVD derivative.

% % Topic: Contribution
To address the aforementioned challenges, we develop an adjoint and RAD-based method to compute the derivative efficiently.
The contribution of the current manuscript is summarized as follows: 
(1) We develop two adjoint-based approaches to compute the derivative of the SVD problem for general imaginary (square and rectangular) matrices.
In the first approach, we leverage the relationship between the EVP of the Gram matrices of the SVD of the target matrix.
In the second approach, we leverage the relationship between the EVP of the target matrix's symmetric embedding and the target matrix's SVD.  
The proposed adjoint techniques can calculate the derivative to machine precision, are straightforward to implement, do not scale with the number of design variables (inputs), and are suitable for use in gradient-based optimization involving SVD.
(2) We develop a RAD-based analytic formula to compute the singular value derivative with respect to a imaginary input matrix, which also works with real-valued matrices in a reduced form.
(3) We propose a generic dot product identity as a tool to derive complex-valued RAD derivative formulae using FAD derivative formulae for imaginary differentiable functions. 

The manuscript is organized as follows.
In~\cref{sec:Governing_equns_and_adjoint_method}, we present the governing equations for the two proposed adjoint method-based approaches for computing the SVD derivative.
Here, the first method relies on forming one of the two Gram-matrices of the main matrix whose SVD derivative is sought.
The second method is based on the symmetric embedding of the target matrix, as mentioned.
In~\cref{sec:der}, we propose these two adjoint-based approaches and a general RAD-based formula for singular value derivative.
In~\cref{sec:Numerical_results}, we compare the results from our proposed adjoint methods and RAD formula for singular value derivative with FD for two randomly selected imaginary matrices, one square and one rectangular.
Next, we show the scalability by testing the singular value derivative on a large dataset sourced from the John Hopkins turbulence database (JHTDB)~\cite{Zaki2013}, which was the transition to turbulence dataset of flow over a flat plate. 
Finally, we present our conclusions in \cref{sec:conc}.
We present a summary of the contributions of researchers in differentiating the SVD and our contributions in \Cref{tab:Summary_table}.

\begin{sidewaystable}[htp]
\centering
\caption{Summary of contributions in literature and current manuscript for SVD derivative computation}\label{tab:Summary_table}
\renewcommand{\arraystretch}{3} % Default value: 1
\begin{tabular}{l l l l l l}
\toprule
Equations & \makecell[l]{Uses only \\ required singular\\ variables\footnote{Requires only $\sigma$, $\mb{u}$ and/or $\mb{v}$ for the derivative calculations, not all the singular variables inside matrices $\mb{U}$, $\mb{\Sigma}$ and/or $\mb{V}$}} & \makecell[l]{Singular \\value \\derivative} & \makecell[l]{Singular \\vector \\derivative} & \makecell[l]{Works for \\complex-valued \\inputs} & Author (s)  \\
\midrule
$\overbar{\mb{A}} = \mb{U} \overbar{\mb{S}} \mb{V}^\intercal$ & \ding{55}& \ding{51} & \ding{55} & \ding{55} & \makecell[l]{\Cref{eq:giles_formula}\\Giles~\cite{Giles2008}}\\
\midrule
\makecell[l]{$\overbar{\mb{A}} = [\mb{U}(\mb{F} \circ [\mb{U}^\intercal \overbar{\mb{U}} - \overbar{\mb{U}}^\intercal \mb{U}])\mb{S}$ \\ $+ (\mb{I}_m - \mb{U}\mb{U}^\intercal)\overbar{\mb{U}}\mb{S}^{-1}]\mb{V}^\intercal$ \\ $+\mb{U}(\mb{I}_k \circ \overbar{\mb{S}}) \mb{V}^\intercal$ \\ $\mb{U}[\mb{S}(\mb{F} \circ [\mb{V}^\intercal \overbar{\mb{V}} - \overbar{\mb{V}}^\intercal \mb{V}])\mb{V}^\intercal$ \\ $+ \mb{S}^{-1} \overbar{\mb{V}}^\intercal (\mb{I}_n - \mb{V}\mb{V}^\intercal)] $ } & \ding{55}& \ding{51} & \ding{51} & \ding{55} & \makecell[l]{\Cref{eq:townsend_formula}\\Townsend~\cite{townsend2016a}} \\
\midrule
\makecell[l]{$\overbar{\mb{A}} = \f{1}{2}(2\mb{U} \overbar{\mb{S}} \mb{V}^* + \mb{U}(\mb{J}+\mb{J}^*)\mb{S}\mb{V}^* $\\
$+ \mb{U}\mb{S}(\mb{K}+\mb{K}^* )\mb{V}^* + \f{1}{2} \mb{U} \mb{S}^{-1}(\mb{L}^*-\mb{L})\mb{V}^*$ \\
$+ 2(\mb{I}-\mb{U}\mb{U}^*)\overbar{\mb{U}}\mb{S}^{-1}\mb{V}^* $\\
$+ 2\mb{U}\mb{S}^{-1} \overbar{\mb{V}}^*(1 - \mb{V} \mb{V}^*))$} & \ding{55} & \ding{51} & \ding{51} & \ding{51} &\makecell[l]{\Cref{eq:wan_and_zhang}\\ Wan and Zhang~\cite{wan2019a}} \\
\midrule
\makecell[l]{$\overbar{\mb{A}} = \mb{U}^\intercal (\mb{G}_2 \mb{V} + \mb{\Lambda}^{-1}\overbar{\mb{V}})$,\\ $\mb{G}_2 = \overbar{\mb{\Lambda}} + 2sym(\mb{G}_1 \circ \mb{E})\mb{\Lambda}$ \\ $- (\mb{\Lambda}^{-1} \overbar{\mb{V}}\mb{V}^\intercal \circ \mb{I})$, \\ $\mb{G}_1 = \overbar{\mb{U}} \mb{U}^\intercal + \mb{\Lambda}^{-1} \overbar{\mb{V}} \mb{V}^\intercal \mb{\Lambda}$} & \ding{55} & \ding{55} & \ding{51} & \ding{55} & \makecell[l]{\Cref{eq:seeger}\\Seeger et al.~\cite{seeger2017a}} \\
\midrule
\makecell[l]{$\mb{M}_g^\intercal \boldsymbol{\psi}_g = (\p g/ \p \mb{w})^\intercal$ \\ $\mb{M}_h^\intercal \boldsymbol{\psi}_h = (\p h/ \p \mb{w})^\intercal$} & \ding{51}& \ding{51} & \ding{51} & \ding{51} & \makecell[l]{\Cref{eq:Adjoint_eqn_for_g,eq:Adjoint_eqn_for_h}} \\
\midrule
\makecell[l]{$\mb{M}_f^\intercal \boldsymbol{\psi}_f = (\p f/ \p \mb{w})^\intercal$} & \ding{51}& \ding{51} & \ding{51} & \ding{51} & \makecell[l]{\Cref{eq:Adjoint_eqn_for_f}} \\
\midrule
\makecell[l]{$\overbar{\mb{A}}_r = (\mb{u}_r \mb{v}_r^\intercal+ \mb{u}_i \mb{v}_i^\intercal)\overbar{\sigma}$ \\ $\overbar{\mb{A}}_i = (-\mb{u}_r \mb{v}_i^\intercal+ \mb{u}_i \mb{v}_r^\intercal)\overbar{\sigma}$}& \ding{51}& \ding{51} & \ding{55} & \ding{51} & \makecell[l]{\Cref{eq:RAD_formula_singular_sens_imaginary}}\\
\bottomrule
\end{tabular}
\end{sidewaystable}

\section{Governing equations}
\label{sec:Governing_equns_and_adjoint_method}
SVD for a general imaginary matrix $\mb{A} \in \mathbb{C}^{m\times n} $ is defined as
\begin{equation}
\mb{A} = \mb{U} \mb{\Sigma} \mb{V}^*,
\label{eq:SVD}
\end{equation}
where $\mb{U} \in \mathbb{C}^{m\times m}$ is the matrix of left singular vectors with each column being a left singular vector, $\mb{V}\in \mathbb{C}^{n\times n}$ is the matrix of right singular vectors with each column being a left singular vector, $\mb{\Sigma} \in \mathbb{R}^{m\times n}$ is a rectangular matrix of singular values on the main diagonal and zeros elsewhere represented by $\sigma_i$ on the main diagonal, arranged in descending order, i.e., $\sigma_1\geq \sigma_2\geq \cdots\geq \sigma_n\geq 0$, $\square^*$ denotes the Hermitian transpose of a imaginary matrix, $m, n\in\mathbb{Z}^+$ are the problem size (without loss of generality, we assume $m>n$).

The SVD and EVP are closely related problems.
There are two ways to compute the singular value via solving an EVP.
The first approach uses the two Gram matrices $\mb{B}=\mb{AA}^* \in \mathbb{C}^{m\times m}$ and $\mb{C}=\mb{A}^*\mb{A} \in \mathbb{C}^{n\times n}$ (refer to~\citet[Page~141]{Schwerdtfeger1961} and to~\citet{Makkonen2023}).
Thus, we name it \textit{Gram matrix method} (GMM).
The second approach is \textit{Symmetric embedding matrix method} (SEMM) (refer to~\citet[Chapter~8.6]{Golub1996} and \citet{Ragnarsson2013}).
We elaborate on these two approaches in the following \cref{sec:GMM,sec:SEMM}.

\subsection{Gram matrix method (GMM)}\label{sec:GMM}
It can be shown that the singular vectors and singular values can be obtained from $\mb{B}$ and $\mb{C}$ (See~\citet[Chapter~8.6]{Golub1996}) 
\begin{equation}
\label{eq:svd_as_evp}
\begin{aligned}
\mb{B} &= \mb{A}\mb{A}^* = \mb{U}\mb{\Sigma}\mb{V}^*\mb{V}\mb{\Sigma}\mb{U}^*=\mb{U}\mb{\Sigma}^2\mb{U}^*,\\
\mb{C} &= \mb{A}^*\mb{A} = \mb{V}\mb{\Sigma}\mb{U}^*\mb{U}\mb{\Sigma}\mb{V}^*=\mb{V}\mb{\Sigma}^2\mb{V}^*.\\
\end{aligned}
\end{equation}
Thus, the left singular vectors, $\mb{U}$, are just the eigenvectors obtained from the eigen-decomposition of $\mb{B}$ and the right singular vectors, $\mb{V}$, are the eigenvectors obtained from the eigen-decomposition of $\mb{C}$. 
The first few singular values are taken from the common eigenvalues between $\mb{B}$ and $\mb{C}$. 
The remaining singular values are taken from the remaining eigenvalues depending on whether $\mb{A}$ is tall or wide.

To apply adjoint method to compute the derivative, we need to construct the residual form of the EVP.
Consider an EVP of a general imaginary square matrix, $\mb{D}\in\mathbb{C}^{m\times m}$, we have
\begin{equation} 
\label{eq:Eigenvalue-decomposition}
\mb{D}\boldsymbol{\phi} = \lambda \boldsymbol{\phi},
\end{equation}
where $\lambda$ and $\boldsymbol{\phi}$ is one eigen-pair.
In general, the eigenvalues and vectors are imaginary. 
It is assumed that the eigenvalues are distinct and hence have multiplicity of one. 
To obtain unique imaginary eigenvectors, we consider the following set of governing equations, in conjunction with \Cref{eq:Eigenvalue-decomposition}.
The residual form is  (see \citet{He2023} for more detail)
\begin{equation}
\label{eq:eigen_gov_equns}
\mb{r}(\mb{w})=
\begin{bmatrix}
\mb{D}\boldsymbol{\phi} - \lambda \boldsymbol{\phi}\\
\boldsymbol{\phi}^* \boldsymbol{\phi} - 1\\
\text{Im}(\boldsymbol{\phi}_k)\\
\end{bmatrix}, \quad 
\mb{w} = 
\begin{bmatrix}
\lambda \\
\boldsymbol{\phi}
\end{bmatrix}, 
\end{equation}
where we further require that $\text{Re}(\boldsymbol{\phi}_k) > 0$, $k=\text{argmax}_j ||\boldsymbol{\phi}_j||_2$, and $\boldsymbol{\phi}_k$ is the $k^{\text{th}}$ element of the vector $\boldsymbol{\phi}$.
This is required because the imaginary eigenvector can always be scaled and rotated in the imaginary plane. 
If we apply stretching and rotation in this plane, we obtain:
\begin{equation}
\mb{D}(\alpha \boldsymbol{\phi} e^{i \theta}) = \lambda(\alpha \boldsymbol{\phi} e^{i \theta}),
\label{eq:eigen_constrain}
\end{equation}
where $\alpha \in \mathbb{R}$ (scaling factor), $\theta \in \mathbb{R}$ (rotation angle). 

The last two equations from \Cref{eq:eigen_gov_equns} therefore constrain the norm of the eigenvector to be one, and constrains the angle of the eigenvector. 

Rewriting the \Cref{eq:eigen_gov_equns} by splitting into their respective real and imaginary parts, we have
\begin{equation}
\label{eq:governing}
\mb{r}(\mb{w}) = 
\begin{bmatrix}
r_{\text{main,r}}\\
r_{\text{main,i}}\\
r_{m}\\
r_{p}
\end{bmatrix}
=
\begin{bmatrix}
\mb{D}_r\boldsymbol{\phi}_r - \mb{D}_i\boldsymbol{\phi}_i - \lambda_r\boldsymbol{\phi}_r + \lambda_i\boldsymbol{\phi}_i\\
\mb{D}_i\boldsymbol{\phi}_r + \mb{D}_r\boldsymbol{\phi}_i - \lambda_i\boldsymbol{\phi}_r - \lambda_r\boldsymbol{\phi}_i\\
\boldsymbol{\phi}_r^\intercal\boldsymbol{\phi}_r + \boldsymbol{\phi}_i^\intercal\boldsymbol{\phi}_i - 1\\
\boldsymbol{e}_k^\intercal \boldsymbol{\phi}_i
\end{bmatrix}, 
\quad
\mb{w} = 
\begin{bmatrix}
\boldsymbol{\phi}_r\\
\boldsymbol{\phi}_i\\
\lambda_r\\
\lambda_i
\end{bmatrix},
\end{equation}
where $\mb{e}_k$ is a vector with all zero elements except that $k^\text{th}$ element equal to one.

For SVD, we can leverage \Cref{eq:svd_as_evp} to convert it to an EVP.
In \Cref{eq:governing}, we can set $\mb{D} = \mb{B}$, then we have a governing equation for the left singular vector and singular value.
In this case, the left singular value becomes the square root of $\lambda$ and the corresponding left singular vector is simply equal to $\boldsymbol{\phi}$ in \Cref{eq:Eigenvalue-decomposition}.
We call this the \textit{Left Gram Matrix Method} (LGMM) since it is the Gram matrix associated with the left singular values.

Alternatively, if we set $\mb{D} = \mb{C}$, we have a governing equation for the right singular vector and singular value.
Again, in this case, the right singular value becomes the square root of $\lambda$ and the corresponding right singular vector is simply equal to $\boldsymbol{\phi}$ in \Cref{eq:Eigenvalue-decomposition}.
We call this the \textit{Right Gram Matrix Method} (RGMM) since it is the Gram matrix associated with the right singular values.

\subsection{Symmetric embedding matrix method}\label{sec:SEMM}
This second approach leverages the following relationship
\begin{equation}
\label{eq:svd_as_evp_2}
\begin{bmatrix}
\mb{O} & \mb{A} \\
\mb{A}^* & \mb{O}
\end{bmatrix}
\begin{bmatrix}
\mb{u}_i\\
\mb{v}_i
\end{bmatrix}
=
\pm\sigma_i
\begin{bmatrix}
\mb{u}_i\\
\mb{v}_i
\end{bmatrix},
\end{equation}
where $\sigma_i, \mb{u}_i$ and $\mb{v}_i$ is one group of singular value and singular vectors (refer to~\citet[Chapter~8.6]{Golub1996} and \citet{Ragnarsson2013}).
Thus, by applying these identities \Cref{eq:svd_as_evp,eq:svd_as_evp_2}, we can obtain the nonlinear governing equations for SVD with additional normalization condition constraints.

In the previous \cref{sec:GMM}, we develop a formulation that relies on the evaluation of $\mb{B}$ or $\mb{C}$ to convert the problem to an EVP where the left or right singular vectors for the corresponding singular value are computed.
The main drawback of that method is that the evaluation of matrix product, $\mb{B}=\mb{A}\mb{A}^*$ or $\mb{C}=\mb{A}^*\mb{A}$, can be expensive when the coefficient matrix is dense.
This motivates us to develop an alternative formulation that is matrix-product free.

We developed the following equation.
By solving this equation, we can get the left, and right singular vectors, and singular value in one shot.
The equation is 
\begin{equation}
\begin{aligned}
\mb{A} \mb{v} &= \sigma \mb{u},\\
\mb{A}^* \mb{u} &= \sigma \mb{v},\\
\mb{u}^*\mb{u}  &= 1,\\
\mathrm{Im}(\mb{u}_k) &= 0.
\end{aligned}
\label{eq:SVD_gov_eqn}
\end{equation}
Notice that by constraining the left singular vector $\mb{u}$, we end up constraining the right singular vector $\mb{v}$.

Starting from \Cref{eq:SVD_gov_eqn}, we can split it into real and imaginary parts and write the residual form as
\begin{equation}
\mb{r}(\mb{w}) = 
\begin{bmatrix}
\mb{A}_r \mb{v}_r - \mb{A}_i \mb{v}_i - \sigma_r \mb{u}_r + \sigma_i \mb{u}_i \\
\mb{A}_r \mb{v}_i + \mb{A}_i \mb{v}_r - \sigma_r \mb{u}_i - \sigma_i \mb{u}_r \\
\mb{A}_r^\intercal \mb{u}_r + \mb{A}_i^\intercal \mb{u}_i - \sigma_r \mb{v}_r + \sigma_i \mb{v}_i \\
\mb{A}_r^\intercal \mb{u}_i - \mb{A}_i^\intercal \mb{u}_r - \sigma_r \mb{v}_i - \sigma_i \mb{v}_r \\
\mb{u}_r^\intercal \mb{u}_r + \mb{u}_i^\intercal \mb{u}_i - 1\\
\mb{e}_k^\intercal \mb{u}_i
\end{bmatrix}, \quad 
\mb{w} = 
\begin{bmatrix}
\mb{u}_r\\
\mb{u}_i\\
\mb{v}_r\\
\mb{v}_i\\
\sigma_r\\
\sigma_i
\end{bmatrix},
\label{eq:residual_form_svd_gov_eqn}
\end{equation}
where $\mb{w}$ is the state vector.
It is well-known that the singular value, $\sigma$, is real. 
However, here we relax it to be a imaginary number to have equal numbers of equations and unknowns.

Now, we prove that the solution of \Cref{eq:SVD_gov_eqn} is indeed a SVD pair.
Pre-multiplying the first equation of \Cref{eq:SVD_gov_eqn} by $\mb{u}^*$, we have
\begin{equation}
\label{eq:1st_premul}
\mb{u}^* \mb{A} \mb{v} = \sigma \mb{u}^*\mb{u} = \sigma.
\end{equation}
Then, we pre-multiply the second equation of \Cref{eq:SVD_gov_eqn} by $\mb{v}^*$, we have
\begin{equation}
\label{eq:2nd_premul}
\mb{v}^* \mb{A}^* \mb{u} = \sigma \mb{v}^*\mb{v} = \sigma ||\mb{v}||_2.
\end{equation}
Then, realizing that the \Cref{eq:1st_premul,eq:2nd_premul} are conjugate transpose, we have

\begin{equation}
\sigma^* = \sigma ||\mb{v}||_2 \Rightarrow \sigma_r = ||\mb{v}||_2 \sigma_r, -\sigma_i=||\mb{v}||_2\sigma_i,
\end{equation}
we have as long as $\sigma = 0$, or $||\mb{v}||_2 = 1$ and $\sigma_i=0$.
Thus, we have implicitly enforced that $\sigma$ is real.

\subsection{Discussion}
The difference between the two formulations GMM and SEMM is that using \Cref{eq:svd_as_evp} method, we can formulate a smaller dimension problem (e.g. if we have $m>n$, we can form $\mb{C}$ which is ${n}\times {n}$).
However, the disadvantage of this method is that it only computes one singular vector (left or right) and singular value, and the other singular vector (right or left) needs to be recovered.
Another disadvantage is that \Cref{eq:svd_as_evp} requires computing matrix-product. 
If the matrix is dense, it can be prohibitively expensive to compute this product.

On the other hand, the conversion using \Cref{eq:svd_as_evp_2} has the advantage that it does not require the matrix product needed in the previous matrix-product formulation.
In addition, the equation includes the singular value and the corresponding left and right singular vectors.
Thus, we do not need to conduct additional computation to search the other singular vector as required by \Cref{eq:svd_as_evp} method.
The disadvantage of this method is that the problem size ($m+n$) is larger than the \Cref{eq:svd_as_evp} method.
We also ensure that there are no repeated singular values else the singular vector space becomes degenerate~\cite{Golub1996,Lin2020}.

\section{Derivative computation}\label{sec:der}
In this section, we propose two general methods to compute the derivative of any function of the singular variables with respect to a general imaginary matrix, $\mb{A}$.
We also propose an RAD formula for just the singular value derivative.
The section is organized as follows:
In \cref{sec:Adjoint_method_summary}, a brief introduction to the adjoint method is given.
\cref{sec:GMM_sens} then gives the adjoint method applied to the governing equations resulting from GMM.
Similarly in \cref{sec:SEMM_sens}, the adjoint method is applied to the governing equations resulting from the SEMM.
Finally, in \cref{sec:RAD}, we propose the RAD formula for the singular value derivative for both imaginary and real-valued cases.

\subsection{Adjoint method}\label{sec:Adjoint_method_summary}

For arbitrary function of interest, $f(\mb{w}, \mb{x})$, we can use the adjoint method to compute its total derivative with respect to $\mb{x}$.
The total derivative is computed using the following formula~\cite{Martins2022}
\begin{equation}
\f{\d f}{\d \mb{x}} = \f{\p f}{\partial \mb{x}} - \boldsymbol{\psi} ^\intercal \f{\partial \mb{r}}{\partial \mb{x}}.
\label{eq:total_der_eqn}
\end{equation}
The adjoint vector, $\boldsymbol{\psi}$, can be found by solving the system of equations resulting from the following adjoint equation
\begin{equation}
\f{\partial \mb{r}}{\partial \mb{w}}^\intercal \boldsymbol{\psi} = \f{\partial f}{\partial \mb{w}}^\intercal,
\label{eq:adjoint_eqn}
\end{equation}
where $\mb{r}=\mb{r}~(\mb{w})$ is the residual vector for the governing equation (\Cref{eq:eigen_gov_equns}), $\mb{w}$ being the state variables vector.
We direct the readers to He et al.~\cite{He2023} for the complete adjoint method to solving the EVP derivative. 
For the summary of the equations and algorithm, see~\Cref{sec:adjoint_evp}.

\hltwo{As a note to the readers, in typical optimization problems, the objective function is $f(\mathbf{w},\mathbf{x})$ and that representation is more general.
In the adjoint approach, we form residual and that is $\mathbf{r}(\mathbf{w},\mathbf{x})$.
The design variables $\mathbf{x}$ affect the matrix $\mathbf{A}$ whose SVD is sought and this in turn affects the objective function $f$.
Since we calculate $\text{d} f / \text{d} \mathbf{A}$ and not $\text{d} f / \text{d} \mathbf{x}$, this is a function $f(\mathbf{w})$ alone.
Therefore, the corresponding residual vector is also a function of $\mb{w}$ as $\mathbf{r}(\mathbf{w})$.}

For the SVD problem, our objective function now becomes:
\begin{equation}
f = f(\mb{u},\mb{v},\sigma, \mb{A}),
\label{eq:svd_obj_func}
\end{equation}
where $f$ can be any function of the singular variables: $\mb{u}$ which is the left singular vector, $\mb{v}$ the right singular vector and these vectors correspond to a chosen singular value (typically the dominant one) $\sigma$.
$\mb{A}$ is also placed in the objective function for the general case if there are any direct dependencies on it, as will be seen in forthcoming sections.

\subsection{Adjoint equation of GMM}\label{sec:GMM_sens}
The computation of the derivative of a general function, $f(\mb{u}, \mb{v}, \sigma, \mb{A})$, can be decomposed into two steps. 
Before we begin with the computation of the derivative of $f(\mb{u},\mb{v},\sigma,\mb{A})$ with respect to $\mb{A}$, consider first the derivatives of $g(\mb{u},\sigma,\mb{A})$ and $h(\mb{v},\sigma,\mb{A})$.
If the function is $g(\mb{u},\sigma,\mb{A})$, it is related to the matrix $\mb{B}$.
If the function is $h(\mb{v},\sigma,\mb{A})$, it is related to the matrix $\mb{C}=\mb{A}^* \mb{A}$ (Refer \cref{sec:GMM}).
We split this function $f$ into its real and imaginary components as $g_r$ and $g_i$ for the function $g$ and $h_r$ and $h_i$ for the function $h$ and therefore we must solve the adjoint for a total of two times for each of $g$ and $h$, once for real part and once for imaginary part of each of $g$ and $h$.
This is done to save computational effort as will be seen in the forthcoming section.
We present this approach as follows.

\subsubsection{Adjoint equation of LGMM}
We found the adjoint equation for the function $g=g(\mb{u},\sigma,\mb{A})\in\mathbb{R}$ to be
\begin{equation}
\mb{M}_{g}^\intercal \boldsymbol{\psi}_{g} = \f{\p g}{\p \mb{w}}^\intercal,
\label{eq:Adjoint_eqn_for_g}
\end{equation}
where $\mb{M}_g$ is
\begin{equation}
\mb{M}_g = 
\begin{bmatrix}
\mb{B}_r - \lambda_r \mb{I} & -\mb{B}_i + \lambda_i \mb{I} & -\mb{u}_{r} & \mb{u}_{i}\\
\mb{B}_i - \lambda_i \mb{I} & \mb{B}_r + \lambda_r \mb{I} & -\mb{u}_{i} & -\mb{u}_{r}\\
2\mb{u}_{r}^\intercal & 2\mb{u}_{i}^\intercal & 0 & 0\\
0 & \mb{e}_{k}^\intercal & 0 & 0
\end{bmatrix},
\label{eq:M_g}
\end{equation}
and $\boldsymbol{\psi}_g$ is
\begin{equation}
\boldsymbol{\psi}_g = 
\begin{bmatrix}
\boldsymbol{\psi}_{\text{$g$,r}} \\
\boldsymbol{\psi}_{\text{$g$,i}} \\
\boldsymbol{\psi}_{\text{m,$g$}}\\
\boldsymbol{\psi}_{\text{p,$g$}}
\end{bmatrix},
\label{eq:psi_g}
\end{equation}
and once the solution to the adjoint vector $\boldsymbol{\psi}$ from either approach is obtained, we can proceed to compute the total derivatives using \Cref{eq:total_der_eqn}. 

We assumed that function $g$ outputs a real number, i.e., $g\in\mathbb{R}$.
If otherwise, function $g$ returns a imaginary number, i.e., $g \in \mathbb{C}$, we can break it into its real and imaginary parts and apply \Cref{eq:Adjoint_eqn_for_g} twice.
We now have two adjoint vectors, namely: $\boldsymbol{\psi}_{g_r}$~and $\boldsymbol{\psi}_{g_i}$~for the function $g$, each with its real and imaginary parts.
The derivative of the residual vector $\mb{r}$ for the function $g$ with respect to the matrix $\mb{B}$ is then  
\begin{equation}
\begin{aligned}
\f{\p \mb{r}_{g_r}}{\p \mb{B}_r}^\intercal \boldsymbol{\psi}_{g_r} &= \boldsymbol{\psi}_{g_r,r} \mb{u}_r^\intercal + \boldsymbol{\psi}_{g_r,i} \mb{u}_i^\intercal,\\
\f{\p \mb{r}_{g_r}}{\p \mb{B}_i}^\intercal \boldsymbol{\psi}_{g_r} &= -\boldsymbol{\psi}_{g_r,r} \mb{u}_i^\intercal + \boldsymbol{\psi}_{g_r,i} \mb{u}_r^\intercal,\\
\f{\p \mb{r}_{g_i}}{\p \mb{B}_r}^\intercal \boldsymbol{\psi}_{g_i} &= \boldsymbol{\psi}_{g_i,r} \mb{u}_r^\intercal + \boldsymbol{\psi}_{g_i,i} \mb{u}_i^\intercal,\\
\f{\p \mb{r}_{g_i}}{\p \mb{B}_i}^\intercal \boldsymbol{\psi}_{g_i} &= -\boldsymbol{\psi}_{g_i,r} \mb{u}_i^\intercal + \boldsymbol{\psi}_{g_i,i} \mb{u}_r^\intercal,
\end{aligned}
\end{equation}
A detailed derivation for this formula can be found in the paper by He et al.~\cite{He2023}.
To proceed to the derivative with respect to the matrix $\mb{A}$, we apply chain rule via the reverse automatic differentiation formula (See~\Cref{sec:Chain_rule_RAD_formula} for derivation).
Thus, to compute the derivatives of the function $f=g$ with respect to matrix $\mb{A}$, we have
\begin{equation}
\begin{aligned}
\f{\d f_r}{\d \mb{A}_r} &= - \biggl(\f{\p \mb{r}_{g_r}}{\p \mb{B}_r}^\intercal \boldsymbol{\psi}_{g_r} \biggl) \mb{A}_r - \biggl( \f{\p \mb{r}_{g_r}}{\p \mb{B}_r}^\intercal \boldsymbol{\psi}_{g_r} \biggl)^\intercal \mb{A}_r - \biggl( \f{\p \mb{r}_{g_r}}{\p \mb{B}_i}^\intercal \boldsymbol{\psi}_{g_r} \biggl)^\intercal \mb{A}_i\\
&+ \biggl( \f{\p \mb{r}_{g_r}}{\p \mb{B}_i}^\intercal \boldsymbol{\psi}_{g_r} \biggl) \mb{A}_i + \f{\p f_r}{\p \mb{A}_r} ,\\
\f{\d f_i}{\d \mb{A}_r} &= -\biggl( \f{\p \mb{r}_{g_i}}{\p \mb{B}_r}^\intercal \boldsymbol{\psi}_{g_i}\biggl) \mb{A}_r - \biggl( \f{\p \mb{r}_{g_i}}{\p \mb{B}_r}^\intercal \boldsymbol{\psi}_{g_i} \biggl)^\intercal \mb{A}_r - \biggl( \f{\p \mb{r}_{g_i}}{\p \mb{B}_i}^\intercal \boldsymbol{\psi}_{g_i} \biggl)^\intercal \mb{A}_i\\
&+ \biggl( \f{\p \mb{r}_{g_i}}{\p \mb{B}_i}^\intercal \boldsymbol{\psi}_{g_i} \biggl) \mb{A}_i + \f{\p f_i}{\p \mb{A}_r} ,\\
\f{\d f_r}{\d \mb{A}_i} &= -\biggl( \f{\p \mb{r}_{g_r}}{\p \mb{B}_r}^\intercal \boldsymbol{\psi}_{g_r} \biggl) \mb{A}_i - \biggl( \f{\p \mb{r}_{g_r}}{\p \mb{B}_r}^\intercal \boldsymbol{\psi}_{g_r} \biggl)^\intercal \mb{A}_i - \biggl(\f{\p \mb{r}_{g_r}}{\p \mb{B}_i}^\intercal \boldsymbol{\psi}_{g_r} \biggl)^\intercal \mb{A}_r\\
&+ \biggl( \f{\p \mb{r}_{g_r}}{\p \mb{B}_i}^\intercal \boldsymbol{\psi}_{g_r} \biggl) \mb{A}_r + \f{\p f_r}{\p \mb{A}_i} ,\\
\f{\d f_i}{\d \mb{A}_i} &= -\biggl( \f{\p \mb{r}_{g_i}}{\p \mb{B}_r}^\intercal \boldsymbol{\psi}_{g_i} \biggl) \mb{A}_i - \biggl(\f{\p \mb{r}_{g_i}}{\p \mb{B}_r}^\intercal \boldsymbol{\psi}_{g_i} \biggl)^\intercal \mb{A}_i - \biggl( \f{\p \mb{r}_{g_i}}{\p \mb{B}_i}^\intercal \boldsymbol{\psi}_{g_i} \biggl)^\intercal \mb{A}_r\\
&+ \biggl( \f{\p \mb{r}_{g_i}}{\p \mb{B}_i}^\intercal \boldsymbol{\psi}_{g_i} \biggl) \mb{A}_r + \f{\p f_i}{\p \mb{A}_i} ,\\
\end{aligned}
\label{eq:get_dfdA_from B}
\end{equation}.
where the last terms in each equation in partial derivative are from the contribution of $\mb{A}$ leading to direct dependencies on it.

As a note to the readers, the derivative expression in~\Cref{eq:get_dfdA_from B} above and many such forthcoming equations involve matrices such as $(\partial \mb{r}_{g_r} / \partial \mb{B}_r)^\top \boldsymbol{\psi}_{g_r}$, which includes complicated tensor-vector products and these can become cumbersome to track and compute. 
To simplify this, we assume the derivative is calculated after flattening the matrices into vectors. 
Once the computation is complete, the resulting vectors are reshaped back to their original matrix dimensions.
The vector flattening and tensor-vector product is elaborated in detail in~\Cref{sec:vectorization}.

Next, consider the objective function $f=f(\mb{u},\mb{v},\sigma,\mb{A})$.
Using the governing equation \Cref{eq:SVD_gov_eqn}, $\mb{u}$ can be expressed as $\mb{A}\mb{v} / \sigma$, then we have that $f(\mb{u}, \mb{v}, \sigma, \mb{A})=f(\mb{u}, \mb{v}(\mb{u}, \sigma, \mb{A}), \sigma, \mb{A})$.
If that is done, then the objective function $f=g(\mb{u},\sigma,\mb{A})$ as shown above.
The same formula in \Cref{eq:get_dfdA_from B} can be used to compute the derivative of $f$ with respect to $\mb{A}$.
This saves the computational effort if $\mb{A}$ is tall.

Since the function $f$ can be any function of $f=f(\mb{u},\mb{v},\sigma,\mb{A})$, there is no general formula that can be devised for that direct contribution of $\mb{A}$.
However, if the only direct contribution of $\mb{A}$ is from changing of $\mb{v}$ to $\mb{A}^* \mb{u}/ \sigma$, then the imaginary partial derivative components of $\p f / \p \mb{A}$ can be written in RAD form shown in~\Cref{sec:RAD_form_for_GMM}.

It is a fact that the last terms of partial derivatives of $\p f / \p \mb{A}$ in \Cref{eq:get_dfdA_from B} can also be obtained by perturbing only $\mb{A}$ in the objective function if present, either using FD or automatic differentiation tools such as JAX~\cite{jax2018github}.
These partial derivative terms in \Cref{eq:get_dfdA_from B} reduce to zero if $f$ is strictly a function of $\mb{u}$ and $\sigma$.

\subsubsection{Adjoint equation of RGMM}
We found the adjoint equation for the function $h=f(\mb{v},\sigma,\mb{A})$ to be
\begin{equation}
\mb{M}_h^\intercal \boldsymbol{\psi}_h = \f{\p h}{\p \mb{w}}^\intercal,
\label{eq:Adjoint_eqn_for_h}
\end{equation}
where $\mb{M}_h$ is
\begin{equation}
\mb{M}_h = 
\begin{bmatrix}
\mb{C}_r - \lambda_r \mb{I} & -\mb{C}_i + \lambda_i \mb{I} & -\mb{v}_{r} & \mb{v}_{i}\\
\mb{C}_i - \lambda_i \mb{I} & \mb{C}_r + \lambda_r \mb{I} & -\mb{v}_{i} & -\mb{v}_{r}\\
2\mb{v}_{r}^\intercal & 2\mb{v}_{i}^\intercal & 0 & 0\\
0 & \mb{e}_{k}^\intercal & 0 & 0
\end{bmatrix},
\label{eq:M_h}
\end{equation}
and $\boldsymbol{\psi}_h$ is
\begin{equation}
\boldsymbol{\psi}_h = 
\begin{bmatrix}
\boldsymbol{\psi}_{\text{$h$,r}} \\
\boldsymbol{\psi}_{\text{$h$,i}} \\
\boldsymbol{\psi}_{\text{m,$h$}}\\
\boldsymbol{\psi}_{\text{p,$h$}}
\end{bmatrix},
\label{eq:psi_h}
\end{equation}
where $h$ can be the real ($h_r$) or imaginary part ($h_i$) of $h$ if it is imaginary.

Note that $\lambda$ which is the eigenvalue is the same in both \Cref{eq:Adjoint_eqn_for_g,eq:Adjoint_eqn_for_h} because both the matrices $\mb{B}$ and $\mb{C}$ will have some common eigenvalues and $\lambda$ is one of them.
($\sigma = \sqrt{\lambda}$).
The system of equations can be solved by a linear solver such as the solver \texttt{numpy.linalg.solve} in \texttt{Python}~\cite{2020NumPy-Array}.
Details on this are shed in He et al~\cite{He2023}.

We now have two adjoint vectors, namely: $\boldsymbol{\psi}_{h_r}$~and $\boldsymbol{\psi}_{h_i}$~for the function $h$, each with its real and imaginary parts.
The derivative of the residual vector $\mb{r}$ for the function $h$ with respect to the matrix $\mb{C}$ is 
\begin{equation}
\begin{aligned}
\f{\p \mb{r}_{h_r}}{\p \mb{C}_r}^\intercal \boldsymbol{\psi}_{h_r} &= \boldsymbol{\psi}_{h_r,r} \mb{v}_r^\intercal + \boldsymbol{\psi}_{h_r,i} \mb{v}_i^\intercal,\\
\f{\p \mb{r}_{h_r}}{\p \mb{C}_i}^\intercal \boldsymbol{\psi}_{h_r} &= -\boldsymbol{\psi}_{h_r,r} \mb{v}_i^\intercal + \boldsymbol{\psi}_{h_r,i} \mb{v}_r^\intercal,\\
\f{\p \mb{r}_{h_i}}{\p \mb{C}_r}^\intercal \boldsymbol{\psi}_{h_i} &= \boldsymbol{\psi}_{h_i,r} \mb{v}_r^\intercal + \boldsymbol{\psi}_{h_i,i} \mb{v}_i^\intercal,\\
\f{\p \mb{r}_{h_i}}{\p \mb{C}_i}^\intercal \boldsymbol{\psi}_{h_i} &= -\boldsymbol{\psi}_{h_i,r} \mb{v}_i^\intercal + \boldsymbol{\psi}_{h_i,i} \mb{v}_r^\intercal.
\end{aligned}
\end{equation}
%Using these derivatives, we compute the derivative of the function $h$, we have
%\begin{equation}
%\begin{aligned}
%\f{\d h_r}{\d \mb{C}_r} = - \f{\p \mb{r}_{h_r}}{\p \mb{C}_r}^\intercal \boldsymbol{\psi}_{h_r},\\
%\f{\d h_r}{\d \mb{C}_i} = - \f{\p \mb{r}_{h_r}}{\p \mb{C}_i}^\intercal \boldsymbol{\psi}_{h_r},\\
%\f{\d h_i}{\d \mb{C}_r} = - \f{\p \mb{r}_{h_i}}{\p \mb{C}_r}^\intercal \boldsymbol{\psi}_{h_i},\\
%\f{\d h_i}{\d \mb{C}_i} = - \f{\p \mb{r}_{h_i}}{\p \mb{C}_i}^\intercal \boldsymbol{\psi}_{h_i}.\\
%\end{aligned}
%\label{eq:get_dfdC_components}
%\end{equation}
%These are the individual derivative components for the function $h$ with respect to $\mb{C}$.
To proceed to the derivative with respect to the matrix $\mb{A}$, we apply chain rule via the reverse automatic differentiation formula (See~\Cref{sec:Chain_rule_RAD_formula} for derivation).
Thus, to compute the derivatives of the function $f=h$ with respect to matrix $\mb{A}$, we have
\begin{equation}
\begin{aligned}
\f{\d f_r}{\d \mb{A}_r} &= -\mb{A}_r \biggl( \f{\p \mb{r}_{h_r}}{\p \mb{C}_r}^\intercal \boldsymbol{\psi}_{h_r} \biggl)^\intercal - \mb{A}_r \biggl( \f{\p \mb{r}_{h_r}}{\p \mb{C}_r}^\intercal \boldsymbol{\psi}_{h_r} \biggl) - \mb{A}_i \biggl( \f{\p \mb{r}_{h_r}}{\p \mb{C}_i}^\intercal \boldsymbol{\psi}_{h_r} \biggl)^\intercal\\
&+ \mb{A}_i \biggl( \f{\p \mb{r}_{h_r}}{\p \mb{C}_i}^\intercal \boldsymbol{\psi}_{h_r} \biggl) + \f{\p f_r}{\p \mb{A}_r},\\
\f{\d f_i}{\d \mb{A}_r} &= -\mb{A}_r \biggl( \f{\p \mb{r}_{h_i}}{\p \mb{C}_r}^\intercal \boldsymbol{\psi}_{h_i} \biggl)^\intercal - \mb{A}_r \biggl( \f{\p \mb{r}_{h_i}}{\p \mb{C}_r}^\intercal \boldsymbol{\psi}_{h_i} \biggl) - \mb{A}_i \biggl( \f{\p \mb{r}_{h_i}}{\p \mb{C}_i}^\intercal \boldsymbol{\psi}_{h_i} \biggl)^\intercal\\
&+ \mb{A}_i \biggl( \f{\p \mb{r}_{h_i}}{\p \mb{C}_i}^\intercal \boldsymbol{\psi}_{h_i} \biggl) + \f{\p f_i}{\p \mb{A}_r},\\
\f{\d f_r}{\d \mb{A}_i} &= -\mb{A}_i \biggl( \f{\p \mb{r}_{h_r}}{\p \mb{C}_r}^\intercal \boldsymbol{\psi}_{h_r} \biggl)^\intercal - \mb{A}_i \biggl( \f{\p \mb{r}_{h_r}}{\p \mb{C}_r}^\intercal \boldsymbol{\psi}_{h_r} \biggl) - \mb{A}_r \biggl( \f{\p \mb{r}_{h_r}}{\p \mb{C}_i}^\intercal \boldsymbol{\psi}_{h_r} \biggl)^\intercal\\
&+ \mb{A}_r \biggl( \f{\p \mb{r}_{h_r}}{\p \mb{C}_i}^\intercal \boldsymbol{\psi}_{h_r} \biggl) + \f{\p f_r}{\p \mb{A}_i},\\
\f{\d f_i}{\d \mb{A}_i} &= -\mb{A}_i \biggl( \f{\p \mb{r}_{h_i}}{\p \mb{C}_r}^\intercal \boldsymbol{\psi}_{h_i} \biggl)^\intercal - \mb{A}_i \biggl( \f{\p \mb{r}_{h_i}}{\p \mb{C}_r}^\intercal \boldsymbol{\psi}_{h_i} \biggl) - \mb{A}_r \biggl( \f{\p \mb{r}_{h_i}}{\p \mb{C}_i}^\intercal \boldsymbol{\psi}_{h_i} \biggl)^\intercal\\
&+ \mb{A}_r \biggl( \f{\p \mb{r}_{h_i}}{\p \mb{C}_i}^\intercal \boldsymbol{\psi}_{h_i} \biggl) + \f{\p f_i}{\p \mb{A}_i},
\end{aligned}
\label{eq:get_dfdA_from_C}
\end{equation}
where the last terms in each equation in partial derivative are from the contribution of $\mb{A}$ leading to direct dependencies on it.

Consider the objective function $f=f(\mb{u},\mb{v},\sigma,\mb{A})$.
Using the governing equation \Cref{eq:SVD_gov_eqn}, $\mb{u}$ can be expressed as $\mb{A}\mb{v} / \sigma$, then we have that $f(\mb{u}, \mb{v}, \sigma, \mb{A})=f(\mb{u}(\mb{v},\sigma,\mb{A}), \mb{v}, \sigma, \mb{A})$.
If that is done, then the objective function $f=h(\mb{v},\sigma,\mb{A})$ as shown above.
The same formula in \Cref{eq:get_dfdA_from_C} can be used to compute the derivative of $f$ with respect to $\mb{A}$.
This saves the computational effort if $\mb{A}$ is wide.

Since the function $f$ can be any function of $f=f(\mb{u},\mb{v},\sigma,\mb{A})$, there is no general formula that can be devised for that direct contribution of $\mb{A}$.
However, if the only direct contribution of $\mb{A}$ is from changing of $\mb{u}$ to $\mb{A} \mb{v}/ \sigma$, then the imaginary partial derivative components of $\p f / \p \mb{A}$ can be written in RAD form shown in~\Cref{sec:RAD_form_for_GMM}.

It is a fact that the last terms of partial derivatives of $\p f / \p \mb{A}$ in \Cref{eq:get_dfdA_from B} can also be obtained by perturbing only $\mb{A}$ in the objective function if present, either using FD or automatic differentiation tools such as JAX~\cite{jax2018github}.
These partial derivative terms in \Cref{eq:get_dfdA_from B} reduce to zero if $f$ is strictly a function of $\mb{u}$ and $\sigma$.

To summarize keeping the minimization of computational effort in mind, we employ \Cref{eq:Adjoint_eqn_for_g,eq:M_g,eq:psi_g,eq:get_dfdA_from B} if the function is expressed in terms of the singular variables $\mb{u}$ and $\sigma$.
If it is expressed in terms of $\mb{v}$ and $\sigma$, then we use \Cref{eq:Adjoint_eqn_for_h,eq:M_h,eq:psi_h,eq:get_dfdA_from_C}.

\subsection{Adjoint equation of SEMM}\label{sec:SEMM_sens}
In this section, we directly apply the adjoint method to the SVD governing equations as listed in the preceding section \cref{sec:Governing_equns_and_adjoint_method}.
We found the adjoint equation for the function, $f$, to be
\begin{equation}
\mb{M}_f^\intercal \boldsymbol{\psi}_f = \f{\p f}{\p \mb{w}}^\intercal,
\label{eq:Adjoint_eqn_for_f}
\end{equation}
\hlthree{where $\mb{M}_f$ and $\boldsymbol{\psi}_f$ are}
\begin{equation}
\mb{M}_f=
\begin{bmatrix}
-\sigma_r \mb{I} & \sigma_i \mb{I} & \mb{A}_r & -\mb{A}_i & -\mb{u}_r & \mb{u}_i\\
-\sigma_i \mb{I} & -\sigma_r \mb{I} & \mb{A}_i & \mb{A}_r & -\mb{u}_i & -\mb{u}_r\\
\mb{A}_r^\intercal & \mb{A}_i^\intercal & -\sigma_r \mb{I} & \sigma_i \mb{I} & -\mb{v}_r & \mb{v}_i\\
-\mb{A}_i^\intercal & \mb{A}_r^\intercal & -\sigma_i \mb{I} & -\sigma_r \mb{I} & -\mb{v}_i & -\mb{v}_r\\
2\mb{u}_r^\intercal & 2\mb{u}_i^\intercal & \mb{0} & \mb{0} & 0 & 0\\
\mb{0} & \mb{e}_k^\intercal & \mb{0} & \mb{0} & 0 & 0
\end{bmatrix},\quad 
\boldsymbol{\psi}_f=
\begin{bmatrix}
\boldsymbol{\psi}_{v_r} \\
\boldsymbol{\psi}_{v_i} \\
\boldsymbol{\psi}_{u_r} \\
\boldsymbol{\psi}_{u_i} \\
\boldsymbol{\psi}_m \\
\boldsymbol{\psi}_p
\end{bmatrix}.\\
\label{eq:SVD_adjoint_eqn}
\end{equation}
Here, $m$ and $p$ represent the magnitude and phase respectively.
The derivatives of $f_r$ and $f_i$ with respect to $\mb{A}_r$ and $\mb{A}_i$ are then  
\begin{equation}
\begin{aligned}
\f{\d f_r}{\d \mb{A}_r} &= -\f{\p \mb{r}}{\p \mb{A}_r}^\intercal \boldsymbol{\psi}_r + \f{\p f_r}{\p \mb{A}_r}\text{~,}\\
\f{\d f_r}{\d \mb{A}_i} &= -\f{\p \mb{r}}{\p \mb{A}_i}^\intercal \boldsymbol{\psi}_r + \f{\p f_r}{\p \mb{A}_i} \text{~,}\\
\f{\d f_i}{\d \mb{A}_r} &= -\f{\p \mb{r}}{\p \mb{A}_r}^\intercal \boldsymbol{\psi}_i + \f{\p f_i}{\p \mb{A}_r} \text{~,}\\
\f{\d f_i}{\d \mb{A}_i} &= -\f{\p \mb{r}}{\p \mb{A}_i}^\intercal \boldsymbol{\psi}_i + \f{\p f_i}{\p \mb{A}_i} \text{~.}
\end{aligned}
\label{eq:dfdA_components_SVD_Adjoint}
\end{equation}
Here, the partial derivatives ${\p f_m}/{\p \mb{A}_n}$ in the RHS of the equations above are for any direct dependencies of $f$ on $\mb{A}$.
The $(\p r / \p \mb{A})^\intercal \boldsymbol{\psi}$ terms in the RHS of the equations above are (derivation in~\Cref{sec:derivation})
\begin{equation}
\begin{aligned}
\f{\p \mb{r}}{\p \mb{A}_r}^\intercal \boldsymbol{\psi} &= \boldsymbol{\psi}_{v_r} \mb{v}_r^\intercal + \boldsymbol{\psi}_{v_i} \mb{v}_i^\intercal + \mb{u}_r \boldsymbol{\psi}_{u_r}^\intercal + \mb{u}_i \boldsymbol{\psi}_{u_i}^\intercal \text{~,}\\
\f{\p \mb{r}}{\p \mb{A}_i}^\intercal \boldsymbol{\psi} &= -\boldsymbol{\psi}_{v_r} \mb{v}_i^\intercal + \boldsymbol{\psi}_{v_i} \mb{v}_r^\intercal + \mb{u}_i \boldsymbol{\psi}_{u_r}^\intercal - \mb{u}_r \boldsymbol{\psi}_{u_i}^\intercal \text{~,}\\
\end{aligned}
\label{eq:get_dfdA_components_direct_SVD}
\end{equation}
where $\boldsymbol{\psi}=\boldsymbol{\psi}_r$ from \Cref{eq:dfdA_components_SVD_Adjoint} if $f = f_r$ in \Cref{eq:SVD_adjoint_eqn}.
Similarly, $\boldsymbol{\psi}=\boldsymbol{\psi}_i$ from \Cref{eq:dfdA_components_SVD_Adjoint} if $f = f_i$ in \Cref{eq:SVD_adjoint_eqn}.

\subsection{RAD method for singular value derivative}\label{sec:RAD}
The technique described so far can compute derivatives of the singular variables at once, in any combination. 
If the derivative of just the dominant singular value is sought, we propose the following RAD based formulae: (Derivation in~\Cref{sec:Singular_Val_Sens_RAD_Derivation}). 
\begin{equation}
\begin{aligned}
\f{\d \sigma}{\d \mb{A}_r} &= \mb{u}_r \mb{v}_r^\intercal + \mb{u}_i \mb{v}_i^\intercal,\\
\f{\d \sigma}{\d \mb{A}_i} &= -\mb{u}_r \mb{v}_i^\intercal + \mb{u}_i \mb{v}_r^\intercal.\\
\end{aligned}
\label{eq:RAD_formula_singular_sens_imaginary}
\end{equation}
For the purely real case if $\mb{A} = \mb{A}_r$ , it is  
\begin{equation}
\f{\d \sigma}{\d \mb{A}} = \mb{u}\mb{v}^\intercal, 
\label{eq:RAD_formula_singular_sens_real}
\end{equation}
where $\mb{u} \in \mathbb{R} $,$\mb{v} \in \mathbb{R}$ and $\mb{A} \in \mathbb{R}$.
Further, the formula in \Cref{eq:RAD_formula_singular_sens_imaginary} reduces to \Cref{eq:RAD_formula_singular_sens_real} if the imaginary parts are set to zero.

\subsection{\hlone{Summary}}\label{sec:summary}

\hlone{In this section we present a summary of the proposed methods in this manuscript.
LGMM, RGMM and SEMM algorithms are shown in} \Cref{LGMM_algorithm}, \Cref{RGMM_algorithm} \hlone{and} \Cref{SEMM_algorithm} \hlone{respectively.}
\begin{algorithm}[H]
    \caption{\hlone{LGMM algorithm}}
    \begin{algorithmic}[1]
    \State Perform SVD on $\mb{A}$ = $\mb{U} \boldsymbol{\Sigma} \mb{V}^*$ and obtain the desired singular variables $(\mb{u}, \sigma, \mb{v})$ % [] TODO HS-: change s-> sigma, S->Sigma, H->* be consistent!! check other symbols.
    \State Construct desired objective function $f$ from the singular variables
    \State Assemble and solve the adjoint equation~(\Cref{eq:Adjoint_eqn_for_g}) to get the adjoint vector $\boldsymbol{\psi}_{g_r}$ from $f_r$ and $\boldsymbol{\psi}_{g_i}$ from $f_i$
    \State Compute the derivatives $\d f_r / \d \mb{A}_r$, $\d f_r / \d \mb{A}_i$, $\d f_i / \d \mb{A}_r$ and $\d f_i / \d \mb{A}_i$ using~\Cref{eq:get_dfdA_from B}
    \end{algorithmic}
    \label{LGMM_algorithm}
\end{algorithm}
\begin{algorithm}[H]
    \caption{\hlone{RGMM algorithm}}
    \begin{algorithmic}[1]
    \State Perform SVD on $\mb{A}$ = $\mb{U} \boldsymbol{\Sigma} \mb{V}^*$ and obtain the desired singular variables $(\mb{u}, \sigma, \mb{v})$
    \State Construct desired objective function $f$ from the singular variables
    \State Assemble and solve the adjoint equation~(\Cref{eq:Adjoint_eqn_for_h}) to get the adjoint vector $\boldsymbol{\psi}_{h_r}$ from $f_r$ and $\boldsymbol{\psi}_{h_i}$ from $f_i$
    \State Compute the derivatives $\d f_r / \d \mb{A}_r$, $\d f_r / \d \mb{A}_i$, $\d f_i / \d \mb{A}_r$ and $\d f_i / \d \mb{A}_i$ using~\Cref{eq:get_dfdA_from_C}
    \end{algorithmic}
    \label{RGMM_algorithm}
\end{algorithm}
\begin{algorithm}[H]
    \caption{\hlone{SEMM algorithm}}
    \begin{algorithmic}[1]
    \State Perform SVD on $\mb{A}$ = $\mb{U} \boldsymbol{\Sigma} \mb{V}^*$ and obtain the desired singular variables $(\mb{u}, \sigma, \mb{v})$
    \State Construct desired objective function $f$ from the singular variables
    \State Assemble and solve the adjoint equation~(\Cref{eq:Adjoint_eqn_for_f}) to get the adjoint vector $\boldsymbol{\psi}_{f_r}$ from $f_r$ and $\boldsymbol{\psi}_{f_i}$ from $f_i$
    \State Compute the derivatives $\d f_r / \d \mb{A}_r$, $\d f_r / \d \mb{A}_i$, $\d f_i / \d \mb{A}_r$ and $\d f_i / \d \mb{A}_i$ using~\Cref{eq:dfdA_components_SVD_Adjoint}
    \end{algorithmic}
    \label{SEMM_algorithm}
\end{algorithm}
\hlone{The RAD equation for only the singular value derivative is given in}~\Cref{eq:Ar_bar_sigma_appendix} \hlone{for the real part of $f$ and in}~\Cref{eq:Ai_bar_sigma_appendix} \hlone{for the imaginary part.}
\hlone{These equations must be used if the objective function $f$ is some function of $\sigma$ alone.
If it is the singular value $\sigma$ itself, then a simpler form is found in}~\Cref{eq:RAD_formula_singular_sens_imaginary} \hlone{for complex-valued problems and is reduced to}~\Cref{eq:RAD_formula_singular_sens_real} \hlone{for real-valued problems.
This finishes our summary of the adjoint and RAD methods proposed in our paper.}

\section{Numerical results}\label{sec:Numerical_results}
In this section, we verify the proposed adjoint-based method via the GMM and SEMM using two complex-valued matrices, one square matrix and one tall rectangular matrix by comparing the values of the obtained derivatives with the values obtained from FD.
The tall rectangular matrix means that it has number of rows greater than the number of columns here~\cite{Trefethen1997a}.
The square matrix is taken and the derivative results from adjoint method for LGMM, RGMM, SEMM and the singular value derivative RAD formula are compared with the results from the FD computation in \cref{sec:Square_matrix_case} for a randomly selected objective function.
This is repeated for the rectangular matrix in \cref{sec:Rect_mat_case}.

Finally, to show scalability, we compute the singular value derivatives of the singular values from the proper orthogonal decomposition (POD) through SVD of a large dataset of transition to turbulence of flow over a flat plate in~\cref{sec:JHTDB}.
The data was sourced from the JHTDB and the POD was performed by the method of snapshots.
\subsection{Square matrix derivative computation}\label{sec:Square_matrix_case}
The randomly selected square matrix is 
\begin{equation}
\mb{A} = \mb{A}_r + i\mb{A}_i = 
\begin{bmatrix}
-1.01 & 0.86 & -31.42\\
3.98 & 0.53 & -7.04\\
3.3 & 8.26 & -3.89
\end{bmatrix}
+ i
\begin{bmatrix}
0.6 & 0.79 & 5.47\\
7.21 & 1.9 & 0.58\\
3.42 & 8.97 & 0.3
\end{bmatrix}.
\label{eq:square_matrix_A}
\end{equation}
It was ensured that no repeated eigenvalues were found in the eigen-decompositions of $\mb{B}=\mb{AA}^*$ and $\mb{C}=\mb{A}^*\mb{A}$.
The dominant singular value and its respective left and right singular vectors are 
\begin{equation}
\begin{aligned}
\sigma_1 &= 33.16357940928816,\\
\mb{u}_{r} + i \mb{u}_{i} &= 
\begin{bmatrix}
0.9572042 \\
0.23641926\\
0.16616206
\end{bmatrix}
+ i
\begin{bmatrix}
0\\
0.01549572\\
0.00401452
\end{bmatrix},\\
\mb{v}_{r} + i \mb{v}_{i} &= 
\begin{bmatrix}
-0.00513976\\
-0.05740808\\
0.99047735
\end{bmatrix}
+ i
\begin{bmatrix}
0.08569246\\
0.09104565\\
0
\end{bmatrix}.\\
\end{aligned}
\label{eq:singular_variables_for_square_mat_A}
\end{equation}
The objective function was chosen as 
\begin{equation}
f = \mb{c}_{\mb{u}} ^ \intercal \mb{u} + \mb{c}_{\mb{v}}^ \intercal \mb{v} + \mathrm{c}_{\sigma}\sigma + \mathrm{c}_{\mb{A}}\text{Tr}(\mb{A}),
\label{eq:func_eg}
\end{equation}
where $\text{Tr}(\cdot)$ is the trace operator and the constants were 
\begin{equation}
\begin{aligned}
\mb{c}_u &= \mb{c}_v =
\begin{bmatrix}
0.16 + 0.78i\\
0.53 + 0.11i\\
0.11 + 0.77i\\
\end{bmatrix},\\
\mathrm{c}_{\sigma} &= \mathrm{c}_\mb{A} = 1,
\end{aligned}
\label{eq:constants}
\end{equation}
where the constants were randomly selected.

In \Cref{sec:GMM_sens}, we showed two ways to compute the derivative.
One is to express $\mb{u}$ in terms of $\mb{v},\sigma,\mb{A}$ which is the LGMM or to express $\mb{v}$ in terms of $\mb{u},\sigma,\mb{A}$ which is the RGMM.
If $\mb{v}$ is expressed in terms of $\mb{u},\sigma,\mb{A}$ following the governing equation in \Cref{eq:SVD_gov_eqn}, the objective function in \Cref{eq:func_eg} becomes 
\begin{equation}
f = \mb{c}_u ^ \intercal \mb{u} + \mb{c}_v ^ \intercal \f{\mb{A}^* \mb{u}}{\sigma} + \mathrm{c}_\sigma\sigma + \mathrm{c}_\mb{A}\text{Tr}(\mb{A}).
\end{equation}
Similarly if $\mb{u}$ is expressed in terms of $\mb{v},\sigma,\mb{A}$, the objective function in \Cref{eq:func_eg} becomes 
\begin{equation}
f = \mb{c}_u ^ \intercal \f{\mb{A} \mb{v}}{\sigma} + \mb{c}_v ^ \intercal \mb{v} + \mathrm{c}_\sigma\sigma + c_\text{A}\text{Tr}(\mb{A}).
\end{equation}
Thus, we have reduced the function from $f(\mb{u},\mb{v},\sigma,\mb{A})$ to either $f(\mb{u},\sigma,\mb{A})$ for LGMM or $f(\mb{v},\sigma,\mb{A})$ for RGMM derivative computation.

For the LGMM, we apply \Cref{eq:Adjoint_eqn_for_g,eq:M_g,eq:psi_g,eq:get_dfdA_from B} to compute the derivative of $f$ with respect to matrix $\mb{A}$.
For the RGMM, we use \Cref{eq:Adjoint_eqn_for_h,eq:M_h,eq:psi_h,eq:get_dfdA_from_C}.
For the \Cref{eq:Adjoint_eqn_for_h} and \Cref{eq:Adjoint_eqn_for_g}, the RHS (Jacobian of the function with respect to the state variables) can be computed analytically, but for the general case, this becomes cumbersome to implement each time the objective function changes.
Thus, the automatic differentiation tool JAX~\cite{jax2018github} was employed for the Jacobian of the function with respect to the state variables. 
This was also done for the direct contribution of matrix $\mb{A}$, as shown in the objective function \Cref{eq:func_eg}, for the partial derivatives in \Cref{eq:get_dfdA_from B,eq:get_dfdA_from_C}.
Thus, the procedure of using JAX~\cite{jax2018github} for the Jacobian and direct contribution of matrix $\mb{A}$ is elaborated in \Cref{sec:Jacobian_and_JAX}.
Results from this derivative computation are shown in \Cref{tab:Adjoint_FD_1}.

Next, the same function and matrix were taken (\Cref{eq:square_matrix_A,eq:func_eg}) and used in the SEMM approach shown in \cref{sec:SEMM_sens}.
The \Cref{eq:Adjoint_eqn_for_f,eq:SVD_adjoint_eqn,eq:dfdA_components_SVD_Adjoint,eq:get_dfdA_components_direct_SVD} were used to compute the derivative.
Once again, JAX~\cite{jax2018github} was used for evaluation of the Jacobian in \Cref{eq:Adjoint_eqn_for_f} and for the partial derivatives of $\p f / \p \mb{A}$ in \Cref{eq:dfdA_components_SVD_Adjoint} (See \Cref{sec:Jacobian_and_JAX}).
Results from this derivative computation are shown in \Cref{tab:Adjoint_FD_1}.

The LGMM, RGMM, and SEMM derivatives were obtained for the same objective function and square matrix and, therefore, matched with each other.
If the objective function $f$ is $f = \sigma$, which is done by setting the other constants to zero and constant $\mathrm{c}_\sigma=1$ in \Cref{eq:func_eg}, then the RAD formula \Cref{eq:RAD_formula_singular_sens_imaginary} can be used to compute the derivative.
Results from this derivative computation are shown in \Cref{tab:RAD_sing_sens_table_Square} for the square matrix $\mb{A}$ (\Cref{eq:square_matrix_A}).
To further establish confidence in the proposed methodology and results, a comparison between the obtained derivative results was made against the results from FD computations (Described in~\Cref{sec:FD_formula}), which matched too, as seen in \Cref{tab:Adjoint_FD_1}.

\subsection{Rectangular matrix derivative computation}\label{sec:Rect_mat_case}
Consider the tall imaginary matrix $\mb{A}$ as
\begin{equation}
\mb{A}_r + i\mb{A}_i =
\begin{bmatrix}
6.3 & 5\\
-5.35 & 0.62\\
-7.49 & -1.6 \\
-0.15 &  0.71\\
\end{bmatrix}
+ i
\begin{bmatrix}
4.49& -9.95\\
-1.23&  7.29\\
6.17& -1.9 \\
-4.89& -3.63\\
\end{bmatrix},
\label{eq:A_rectangular}
\end{equation}
and the chosen singular variables are
\begin{equation}
\begin{aligned}
\sigma &= 17.275386033399094,\\
\mb{u}_{r} + i \mb{u}_{i} &= 
\begin{bmatrix}
-0.64373688\\
0.51599116\\
0.2346775\\
-0.13537246\\
\end{bmatrix}
+ i
\begin{bmatrix}
-0.41836418\\
0.05007437\\
-0.20205518\\
0.16611561\\
\end{bmatrix},\\
\mb{v}_{r} + i \mb{v}_{i} &= 
\begin{bmatrix}
-0.72661509\\
0.05431442\\
\end{bmatrix}
+ i
\begin{bmatrix}
0\\
-0.68489449\\
\end{bmatrix}.
\end{aligned}
\label{eq:singular_vars_rect_mat}
\end{equation}
The objective function was taken to be same as in \Cref{eq:func_eg} and the constants were the same as well (\Cref{eq:constants}), except for $\mb{c}_u$ and $\mb{c}_v$ which were chosen to be 
\begin{equation}
\mb{c}_u = 
\begin{bmatrix}
0.12+0.67i\\
0.56+3.67i\\
0.46+2.96i\\
2.89+1.48i
\end{bmatrix},
\text{~} \mb{c}_v = 
\begin{bmatrix}
7.12+0.97i\\
0.26+6.47i\\
\end{bmatrix}.
\end{equation}
Once again, these constants were randomly selected.
The constants  $\mb{c}_u$ and $\mb{c}_v$ for the rectangular matrix were chosen differently because they must adhere to the shape of the vectors $\mb{u}$ and $\mb{v}$.

The same process as in \cref{sec:Square_matrix_case} was followed for the derivative computation of each of the routes: LGMM, RGMM and SEMM, and the results were compared with FD.
The results are tabulated in \Cref{tab:Adjoint_FD_2}.

If the objective function $f$ is $f = \sigma$, which is done by setting the other constants to zero and constant $\mathrm{c}_\sigma=1$ in \Cref{eq:func_eg}, then the RAD formula \Cref{eq:RAD_formula_singular_sens_imaginary} can be employed.
Results from this derivative computation are shown in \Cref{tab:RAD_sing_sens_table_rect} for the rectangular matrix $\mb{A}$ (\Cref{eq:A_rectangular}).
\begin{sidewaystable}
\centering
\caption{Verification of the adjoint derivatives with FD for square matrix case}
\begin{tabular}{llllll}
\toprule
Type
&Index
&\multicolumn{1}{c}{\makecell{Adjoint \\ LGMM}}
&\multicolumn{1}{c}{\makecell{Adjoint \\ RGMM}}
&\multicolumn{1}{c}{\makecell{Adjoint \\ SEMM}}
&\multicolumn{1}{c}{FD} \\
\midrule
\multicolumn{1}{c}{${\d f_r}/{\d \mb{A}_r}$} & \multicolumn{1}{c}{$(1, 1)$} & \multicolumn{1}{r}{$1.006352961803713$} & \multicolumn{1}{r}{$1.006352961803713$} & \multicolumn{1}{r}{$1.006352961803713$} &\multicolumn{1}{r}{$1.006352\underline{068344540}$}\\
\multicolumn{1}{c}{${\d f_r}/{\d \mb{A}_r}$} & \multicolumn{1}{c}{$(1, 2)$} & \multicolumn{1}{r}{$0.043276271008604$} & \multicolumn{1}{r}{$0.043276271008604$} & \multicolumn{1}{r}{$0.043276271008604$} &\multicolumn{1}{r}{$0.04327\underline{5413474930}$}\\
\multicolumn{1}{c}{${\d f_r}/{\d \mb{A}_r}$} & \multicolumn{1}{c}{$(1, 3)$} & \multicolumn{1}{r}{$-0.936930641525170$} & \multicolumn{1}{r}{$-0.936930641525170$} & \multicolumn{1}{r}{$-0.936930641525170$} &\multicolumn{1}{r}{$-0.93693\underline{0476314046}$}\\
\midrule
\multicolumn{1}{c}{${\d f_r}/{\d \mb{A}_i}$} & \multicolumn{1}{c}{$(1, 1)$} & \multicolumn{1}{r}{$0.063695888970744$} & \multicolumn{1}{r}{$0.063695888970744$} & \multicolumn{1}{r}{$0.063695888970744$} &\multicolumn{1}{r}{$0.06369\underline{6809604608}$}\\
\multicolumn{1}{c}{${\d f_r}/{\d \mb{A}_i}$} & \multicolumn{1}{c}{$(1, 2)$} & \multicolumn{1}{r}{$0.082766443637231$} & \multicolumn{1}{r}{$0.082766443637231$} & \multicolumn{1}{r}{$0.082766443637231$} &\multicolumn{1}{r}{$0.08276\underline{7300568776}$}\\
\multicolumn{1}{c}{${\d f_r}/{\d \mb{A}_i}$} & \multicolumn{1}{c}{$(1, 3)$} & \multicolumn{1}{r}{$0.156944460318835$} & \multicolumn{1}{r}{$0.156944460318835$} & \multicolumn{1}{r}{$0.156944460318835$} &\multicolumn{1}{r}{$0.15694\underline{6299512128}$}\\
\midrule
\multicolumn{1}{c}{${\d f_i}/{\d \mb{A}_r}$} & \multicolumn{1}{c}{$(1, 1)$} & \multicolumn{1}{r}{$-0.017846334274906$} & \multicolumn{1}{r}{$-0.017846334274906$} & \multicolumn{1}{r}{$-0.017846334274906$} &\multicolumn{1}{r}{$-0.017846\underline{180533354}$}\\
\multicolumn{1}{c}{${\d f_i}/{\d \mb{A}_r}$} & \multicolumn{1}{c}{$(1, 2)$} & \multicolumn{1}{r}{$0.002833157022766$} & \multicolumn{1}{r}{$0.002833157022766$} & \multicolumn{1}{r}{$0.002833157022766$} &\multicolumn{1}{r}{$0.002833\underline{298928806}$}\\
\multicolumn{1}{c}{${\d f_i}/{\d \mb{A}_r}$} & \multicolumn{1}{c}{$(1, 3)$} & \multicolumn{1}{r}{$0.006354411136774$} & \multicolumn{1}{r}{$0.006354411136774$} & \multicolumn{1}{r}{$0.006354411136774$} &\multicolumn{1}{r}{$0.006354\underline{630599503}$}\\
\midrule
\multicolumn{1}{c}{${\d f_i}/{\d \mb{A}_i}$} & \multicolumn{1}{c}{$(1, 1)$} & \multicolumn{1}{r}{$1.006719107202561$} & \multicolumn{1}{r}{$1.006719107202561$} & \multicolumn{1}{r}{$1.00671910 7202561$} &\multicolumn{1}{r}{$1.00671910\underline{0082876}$}\\
\multicolumn{1}{c}{${\d f_i}/{\d \mb{A}_i}$} & \multicolumn{1}{c}{$(1, 2)$} & \multicolumn{1}{r}{$0.019954391307143$} & \multicolumn{1}{r}{$0.019954391307143$} & \multicolumn{1}{r}{$0.019954391307143$} &\multicolumn{1}{r}{$0.019954\underline{295993330}$}\\
\multicolumn{1}{c}{${\d f_i}/{\d \mb{A}_i}$} & \multicolumn{1}{c}{$(1, 3)$} & \multicolumn{1}{r}{$0.003910503966226$} & \multicolumn{1}{r}{$0.003910503966226$} & \multicolumn{1}{r}{$0.0039105 03966226$} &\multicolumn{1}{r}{$0.0039105\underline{11914285}$}\\
\bottomrule
\end{tabular}\label{tab:Adjoint_FD_1}
\end{sidewaystable}

\begin{sidewaystable}
\centering
\caption{Verification of the adjoint derivatives with FD for rectangular matrix case}
\begin{tabular}{llllll}
\toprule
Type
&Index
&\multicolumn{1}{c}{\makecell{Adjoint \\ LGMM}}
&\multicolumn{1}{c}{\makecell{Adjoint \\ RGMM}}
&\multicolumn{1}{c}{\makecell{Adjoint \\ SEMM}}
&\multicolumn{1}{c}{FD} \\
\midrule
\multicolumn{1}{c}{${\d f_r}/{\d \mb{A}_r}$} & \multicolumn{1}{c}{$(1, 1)$} & \multicolumn{1}{r}{$1.846102900714162$} & \multicolumn{1}{r}{$1.846102900714162$} & \multicolumn{1}{r}{$1.846102900714162$} &\multicolumn{1}{r}{$1.84610\underline{1064018058}$}\\
\multicolumn{1}{c}{${\d f_r}/{\d \mb{A}_r}$} & \multicolumn{1}{c}{$(1, 2)$} & \multicolumn{1}{r}{$-0.006821647620363$} & \multicolumn{1}{r}{$-0.006821647620363$} & \multicolumn{1}{r}{$-0.006821647620363$} &\multicolumn{1}{r}{$-0.00682\underline{2084230862}$}\\
\midrule
\multicolumn{1}{c}{${\d f_r}/{\d \mb{A}_i}$} & \multicolumn{1}{c}{$(1, 1)$} & \multicolumn{1}{r}{$0.354647919899091$} & \multicolumn{1}{r}{$0.354647919899091$} & \multicolumn{1}{r}{$0.354647919899091$} &\multicolumn{1}{r}{$0.354647\underline{895051130}$}\\
\multicolumn{1}{c}{${\d f_r}/{\d \mb{A}_i}$} & \multicolumn{1}{c}{$(1, 2)$} & \multicolumn{1}{r}{$-0.124582311966832$} & \multicolumn{1}{r}{$-0.124582311966832$} & \multicolumn{1}{r}{$-0.124582311966832$} &\multicolumn{1}{r}{$-0.12458\underline{1340799068}$}\\
\midrule
\multicolumn{1}{c}{${\d f_i}/{\d \mb{A}_r}$} & \multicolumn{1}{c}{$(1, 1)$} & \multicolumn{1}{r}{$0.227870193134751$} & \multicolumn{1}{r}{$0.227870193134751$} & \multicolumn{1}{r}{$0.227870193134751$} &\multicolumn{1}{r}{$0.2278701\underline{47639919}$}\\
\multicolumn{1}{c}{${\d f_i}/{\d \mb{A}_r}$} & \multicolumn{1}{c}{$(1, 2)$} & \multicolumn{1}{r}{$0.353104252473483$} & \multicolumn{1}{r}{$0.353104252473483$} & \multicolumn{1}{r}{$0.353104252473483$} &\multicolumn{1}{r}{$0.35310\underline{3555283951}$}\\
\midrule
\multicolumn{1}{c}{${\d f_i}/{\d \mb{A}_i}$} & \multicolumn{1}{c}{$(1, 1)$} & \multicolumn{1}{r}{$0.780271281223158$} & \multicolumn{1}{r}{$0.780271281223158$} & \multicolumn{1}{r}{$0.780271281223158$} &\multicolumn{1}{r}{$0.78027\underline{3927247777}$}\\
\multicolumn{1}{c}{${\d f_i}/{\d \mb{A}_i}$} & \multicolumn{1}{c}{$(1, 2)$} & \multicolumn{1}{r}{$0.113513089065063$} & \multicolumn{1}{r}{$0.113513089065063$} & \multicolumn{1}{r}{$0.113513089065063$} &\multicolumn{1}{r}{$0.11351\underline{2610866451}$}\\
\bottomrule
\end{tabular}
\label{tab:Adjoint_FD_2}
\end{sidewaystable}

\begin{table}[H]
\centering
\caption{Verification of the RAD formula for singular value derivative with FD for the square matrix}
\begin{tabular}{llll}
\toprule
Type
&Index
&\multicolumn{1}{c}{RAD}
&\multicolumn{1}{c}{FD} \\
\midrule
\multicolumn{1}{c}{${\d \sigma}/{\d \mb{A}_r}$} & \multicolumn{1}{c}{$(1, 1)$} & \multicolumn{1}{r}{$0.018703061899253$}  & \multicolumn{1}{r}{$0.01870\underline{2181137087}$}\\
\multicolumn{1}{c}{${\d \sigma}/{\d \mb{A}_r}$} & \multicolumn{1}{c}{$(1, 2)$} & \multicolumn{1}{r}{$0.068881276214858$}  & \multicolumn{1}{r}{$0.06888\underline{0360970525}$}\\
\multicolumn{1}{c}{${\d \sigma}/{\d \mb{A}_r}$} & \multicolumn{1}{c}{$(1, 3)$} & \multicolumn{1}{r}{$-0.934470093986586$}  & \multicolumn{1}{r}{$-0.9344700\underline{51561220}$}\\
\midrule
\multicolumn{1}{c}{${\d \sigma}/{\d \mb{A}_i}$} & \multicolumn{1}{c}{$(1, 1)$} & \multicolumn{1}{r}{$0.080015153716675$}  & \multicolumn{1}{r}{$0.08001\underline{6050674203}$}\\
\multicolumn{1}{c}{${\d \sigma}/{\d \mb{A}_i}$} & \multicolumn{1}{c}{$(1, 2)$} & \multicolumn{1}{r}{$0.076615998520789$}  
& \multicolumn{1}{r}{$0.07661\underline{6849753464}$}\\
\multicolumn{1}{c}{${\d \sigma}/{\d \mb{A}_i}$} & \multicolumn{1}{c}{$(1, 3)$} & \multicolumn{1}{r}{$0.160118791389606$}  
& \multicolumn{1}{r}{$0.1601\underline{20613657000}$}\\
\bottomrule
\end{tabular}
\label{tab:RAD_sing_sens_table_Square}
\end{table}

\begin{table}[H]
\centering
\caption{Verification of the RAD formula for singular value derivative with FD for the rectangular matrix}
\begin{tabular}{llll}
\toprule
Type
&Index
&\multicolumn{1}{c}{RAD}
&\multicolumn{1}{c}{FD} \\
\midrule
\multicolumn{1}{c}{${\d \sigma}/{\d \mb{A}_r}$} & \multicolumn{1}{c}{$(1, 1)$} & \multicolumn{1}{r}{$0.467749108787955$}  
& \multicolumn{1}{r}{$0.46774\underline{8947130531}$}\\
\multicolumn{1}{c}{${\d \sigma}/{\d \mb{A}_r}$} & \multicolumn{1}{c}{$(1, 2)$} & \multicolumn{1}{r}{$0.251572392322310$}  
& \multicolumn{1}{r}{$0.25157\underline{1147913410}$}\\
\midrule
\multicolumn{1}{c}{${\d \sigma}/{\d \mb{A}_i}$} & \multicolumn{1}{c}{$(1, 1)$} & \multicolumn{1}{r}{$0.303989439817427$}  
& \multicolumn{1}{r}{$0.303989\underline{750705114}$}\\
\multicolumn{1}{c}{${\d \sigma}/{\d \mb{A}_i}$} & \multicolumn{1}{c}{$(1, 2)$} & \multicolumn{1}{r}{$-0.463615299320613$}  
& \multicolumn{1}{r}{$-0.463615\underline{023704733}$}\\
\bottomrule
\end{tabular}
\label{tab:RAD_sing_sens_table_rect}
\end{table}

\subsection{\hlo{Scalability of the proposed methodology}}\label{sec:scalability}
\hlo{In this section, the computational speed and scalability study of the proposed methodology was conducted by measuring the mean wall time taken for real-valued LGMM, RGMM, SEMM, FD and RAD, for multiple square matrices of increasing matrix size.
Complex-valued would yield similar results.
However, the RAD-formula from~Townsend~\cite{townsend2016a} in}~\Cref{eq:townsend_formula} \hlo{is for purely real-valued matrices and therefore we elected to go with a purely real-valued comparison.
The adjoint equations were solved using LAPACK as part of the linalg library of NumPy in Python~\cite{2020NumPy-Array}.}

\hlo{The matrices were ensured to be full rank matrices and randomized for each value of the $N \times N$ matrix.
In all of the derivative approaches, the computational wall time measurement was done carefully.}

\hlo{For the real-valued part, only the real-valued parts of LGMM, RGMM and SEMM were compared against that of the real-valued RAD formula for SVD derivative formula proposed by} Townsend~\cite{townsend2016a} \Cref{eq:Adjoint_eqn_for_g,eq:Adjoint_eqn_for_g,eq:Adjoint_eqn_for_f,eq:townsend_formula}.

\hlo{It can be seen that for just the real-valued part, the RAD formula is relatively expensive in computational effort compared to the adjoint.
Further, efficient preconditioning of the solver for the adjoint equation could improve the speed of the computation (See Section~4.3 from He et al.~\cite{He2023} for more details).}
\begin{figure}[H]
    \centering
    \includegraphics[width=1\textwidth]{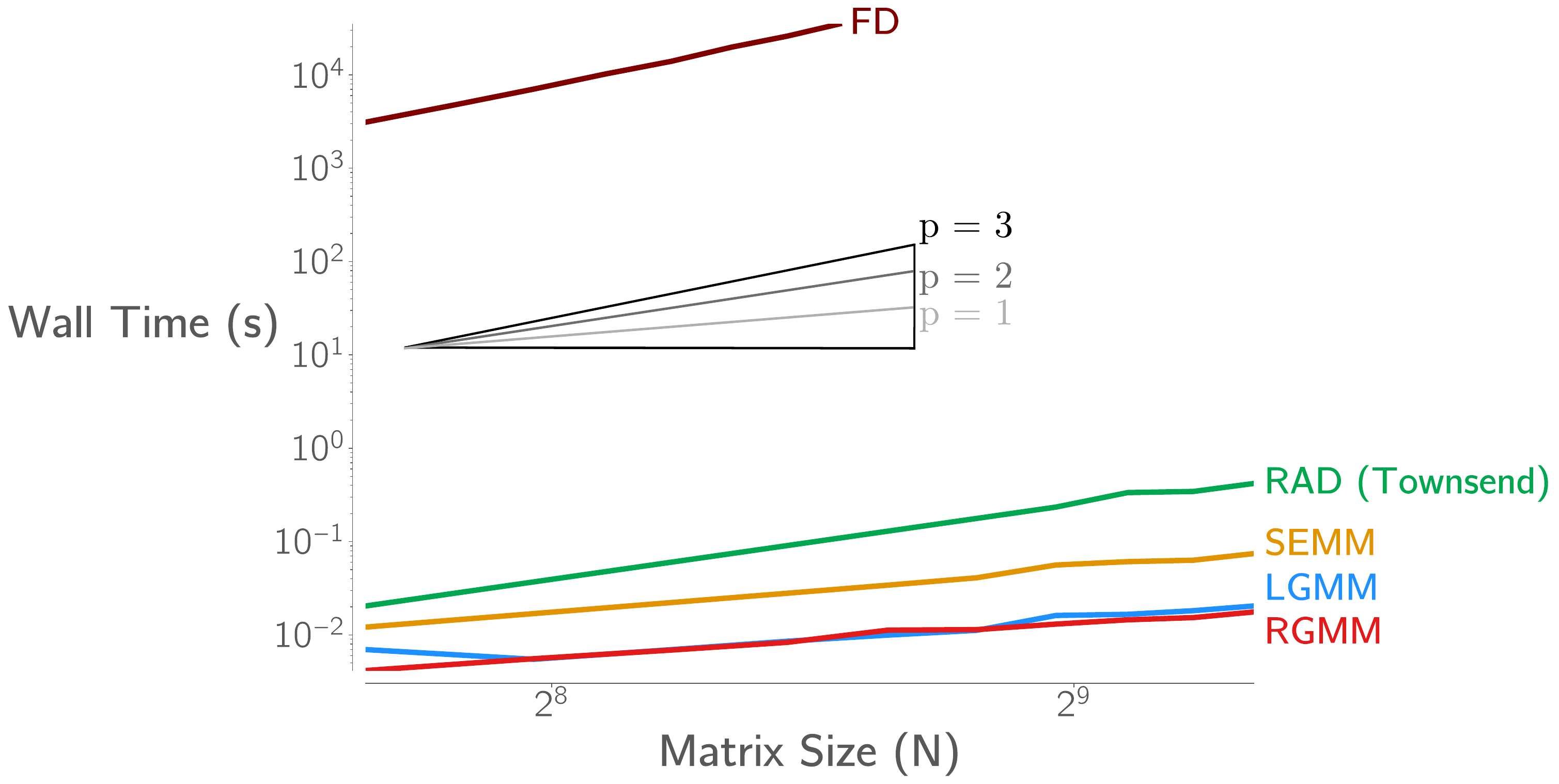}
    \caption{\hlo{Computational efficiency comparison of real-valued proposed methods, FD and RAD.
    $p$ represents the index for slope as a reference and $p=1$ corresponds to our proposed methods, $p=2$ is for RAD and $p=3$ is for FD in the plot.
    The smaller slope and lower wall time, yields better scalability.}}
    \label{fig:RAD_comparison}
\end{figure}

\hlo{For the RAD-adjoint comparison study performed, the imaginary part of the matrix is just a zero-matrix, but is inevitably present in the adjoint equation, for instance, as} $\mb{A}_i$ in~\Cref{eq:SVD_adjoint_eqn}.
\hltwo{This simply increases the sparsity of the adjoint coefficients matrix $\mb{M}_f$ from}~\Cref{eq:Adjoint_eqn_for_f}.
\hltwo{A sparse solver would therfore be even faster. %[x] TODO HS-: bold for matrix and define it.
Alternately a purely real-valued formulation for the adjoint method would result in lesser computational time.
Even with this, the adjoint method performs better.
Further, it is noted that for the comparison cases of proposed methods with FD, the SVD computation time was immaterial as as the same algorithm as in SVD solver NumPy~\cite{2020NumPy-Array} was used.}

% [x] TODO HS-: ref 11 looks strange fix.

\hlo{However, for the case of RAD comparison as shown in} \cref{fig:RAD_comparison}, \hlo{the RAD method requires all the singular variables whereas our methods require only the singular variables under question.
Thus, the computation time for the full SVD through NumPy was added to the computational time measured for RAD, while for our methods, a real-valued non-linear iterative solver was made to solve}~\Cref{eq:SVD_gov_eqn}.
\hlo{The computational time to run this for only the singular variables under consideration was measured and added to the times measured for the derivative computations for our methods.}

\hlo{Finally, it is evident from}~\cref{fig:RAD_comparison} \hlo{that the FD scales poorly with the number of inputs and during the computational study of comparsion of these methods, the matrix size could not be increased beyond the final size attained in}~\cref{fig:RAD_comparison}\hlo{ as the time taken was of the order $10^5$ seconds, and thus not feasible, making our case of scalability issue with FD.}

\subsection{Large dataset derivative computation}\label{sec:JHTDB}
In this section, we compute the derivative of the first six singular values from the proper orthogonal decomposition of the fluid flow over a flat plate. 
This flow is transitioning from laminar to turbulent regime and is sourced from JHTDB~\cite{Zaki2013}. 
We use the RAD based singular value derivative for the real-valued case~\Cref{eq:RAD_formula_singular_sens_real} on the JHTDB dataset.

\subsubsection{Dataset description}\label{sec:Dataset_descr}
The laminar--turbulent transitional flow dataset from JHTB was used in this study~\cite{Zaki2013}.
This data-set was produced from direct numerical simulation of incompressible flow of fluid over a flat plate in a developing boundary layer.
The flat plate was given a leading edge and the Reynolds number $Re=800$, where $Re = U_\infty L/\nu$.
Here, $U_\infty=1$ is the non-dimensional free stream velocity, $L=1$ is the length scale set to half-plate thickness and $\nu=1.24\times10^{-3}$ is the kinematic viscosity.

The dataset is a rectangular domain in 3D with set of grid points in cartesian coordinates $(x,y,z)$, each having the state variables of $x$-component of velocity $u$, $y$-component of velocity $v$, $z$-component of velocity $w$ and the pressure $P$. 
All values in the database are non-dimensionalized.
The database domain is from $[30.2185L,1000L]$ in the $x$ direction, $[0,26.48L]$ in the $y$ direction and $[0,240L]$ in the $z$ direction.
The cutout was taken from $[100L,600L]$ in the $x$ direction, $[0,10L]$ in the $y$ direction and $[80L,120L]$ in the $z$ direction, which resulted in a domain size of $(500L \times 10L \times 40L)$.
The number of grid points in the cutout was $n_x\times n_y\times n_z = 1712 \times 85 \times 342$.
Here $n_x, n_y$ and $n_z$ are the number of grid points in the $x$, $y$ and $z$ directions, respectively.
Thus, there are $m_0=1712 \times 85 \times 342=49.76784\times 10^6$ spatial points.
The cutout time stored is $\bar{t} \in [1000,1175]$ where $\bar{t}=L/U_\infty$ is the non-dimensional time.
Thus, there were $75$ time-steps stored, with each step having a physical time of $\delta \bar{t}=0.25$ giving a total physical time, $\bar{t}$, equal to $75\times0.25=18.75$.
The last 75 time steps were chosen after studying the flow field carefully and ensuring a statistically stationary state, where by stationary we mean the transition location.
Further details on the grid and velocity data, the turbulence characteristics and other flow related data can be found at JHTDB~\cite{Zaki2013}.

\subsubsection{POD and its relationship with SVD}\label{sec:POD_description}
POD of the obtained data cutout was performed on the velocity data.
Thus, there are $m=n_x \times n_y \times n_z \times 3 = m_0 \times 3 = 149.30352\times 10^{6}$ states for each time step.
In POD, the time-series data is arranged into a snapshot matrix $\mb{X} \in\mathbb{R}^{m\times n}$ where $m$ is the number of states arranged in a vectorized form for each time step and $n$ equals to the number of time-steps or snapshots of the data. 
SVD on this snapshot matrix is performed to obtain the modes, their energies and their temporal coefficients.
In the current study, since velocity data was used for the POD, the matrix $\mb{X}$ is defined as 
\begin{equation}
\mb{X} = 
\begin{bmatrix}
u_{1,1} & u_{1,2} & \cdots & u_{1,n}\\
u_{2,1} & u_{2,2} & \cdots & u_{2,n}\\
\vdots & \vdots & \vdots &\vdots\\
u_{m_0,1} & u_{m_0,2} & \cdots & u_{m_0,n}\\
\makebox[0pt][l]{\hspace{-0.5cm}\rule{0.331\linewidth}{0.4pt}}\\
v_{1,1} & v_{1,2} & \cdots & v_{1,n}\\
v_{2,1} & v_{2,2} & \cdots & v_{2,n}\\
\vdots & \vdots & \vdots &\vdots\\
v_{m_0,1} & v_{m_0,2} & \cdots & v_{m_0,n}\\
\makebox[0pt][l]{\hspace{-0.5cm}\rule{0.331\linewidth}{0.4pt}}\\
w_{1,1} & w_{1,2} & \cdots & w_{1,n}\\
w_{2,1} & w_{2,2} & \cdots & w_{2,n}\\
\vdots & \vdots & \vdots &\vdots\\
w_{m_0,1} & w_{m_0,2} & \cdots & w_{m_0,n}\\
\end{bmatrix}
,
\label{eq:snapshot_matrix_representation}
\end{equation}
where $u_{i,k}$ is the $x$-directional velocity at the $i^\text{th}$ grid point and $k^\text{th}$ time-step; similar for $v_{i,k}$ and $w_{i,k}$.
The number of time-steps is $n=20$ as mentioned earlier. 
Instead of conducting POD on the original dataset, it is recommended to instead conduct the POD with respect to the perturbed state~\cite{Taira2020}.
The perturbed snapshot matrix is defined as
\begin{equation}
\mb{X}' = \mb{X} - \f{1}{n} \mb{X} \mb{1}_{n} \mb{1}_{n}^\intercal,
\end{equation}
where $\mb{X}'\in\mathbb{R}^{m\times n}$ is the perturbed snapshots matrix, and $\mb{1}_{n}\in \mathbb{R}^{n}$ is a vector of ones with dimension $n$.
SVD is then performed on this matrix $\mb{X}'$ which yields
\begin{equation}
\mb{X}' = \mb{\Phi} \mb{\Sigma} \boldsymbol{\Psi}^\intercal,
\end{equation}
where $\boldsymbol{\Phi} \in\mathbb{R}^{m \times m}$, $\mb{\Sigma} \in\mathbb{R}^{m \times n}$, and $\boldsymbol{\Psi} \in\mathbb{R}^{n \times n}$.
Columns of $\mb{\Phi}$ contain the POD modes, $\mb{\Sigma}$ is the rectangular matrix with singular values on the main diagonal and zeros elsewhere, and each of these values contains the energy associated with the corresponding mode.
The POD was performed through the method of snapshots (See~\Cref{sec:mthod_of_snapshots}) that gave the singular values and vectors as it is an efficient way to perform POD~\cite{Taira2020}.

\subsubsection{POD derivative results}\label{sec:POD_diff}
\begin{figure}[H]
\centering
\includegraphics[width=0.8\textwidth]{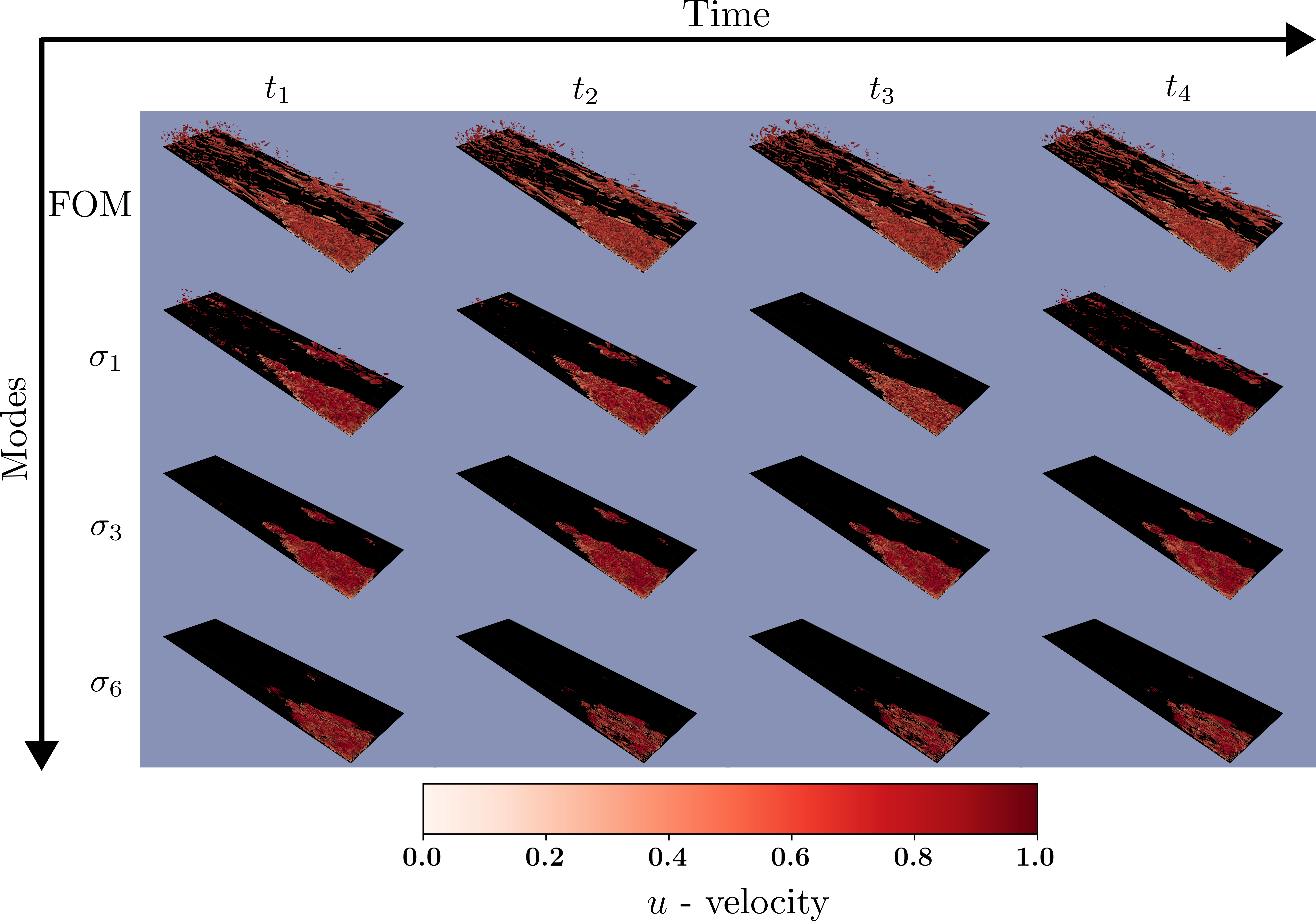}
\caption{POD modes of the JHTDB transitional boundary layer dataset. Vortical structures are \(Q\) iso-surfaces (\(Q\) = $0.001$)}
\label{fig:modes}
\end{figure}

The results of the POD computations along with the full order model (FOM) are shown in~\cref{fig:modes}.
The first, third and sixth modes are shown for sake of compactness and representation, although all of the modes were captured in the computations.
The velocity data was visualized using the \(Q\)-criterion with \(Q\)-2 iso-surfaces of $0.001$ in magnitude and colored by the $u$-velocity magnitude of the FOM flow field \cite{hunt1987}.
The \(Q\)-criterion identifies vortices as regions where the rotation rate dominates over the strain rate. It is defined as
\begin{equation}
Q = \frac{1}{2} \left( ||\mathbf{\Omega}||^2 - ||\mathbf{S}||^2 \right),
\end{equation}
where \(Q\) is the \(Q\)-criterion value, \(\mathbf{\Omega}\) is the antisymmetric rotation tensor, $\mathbf{S}$ is the symmetric strain-rate tensor, $||\mathbf{\Omega}||^2$ is the squared norm of the rotation tensor, and $||\mathbf{S}||^2$ is the squared norm of the strain-rate tensor. These tensors are defined as
\begin{equation}
\mathbf{S} = \frac{1}{2} \left( \nabla \mathbf{u} + \nabla \mathbf{u}^\intercal \right),
\end{equation}
where $\mb{u}$ is the velocity vector containing the $x,y$ and $z$ component of velocities in a vectorized form and $\nabla \mb{u}$ is the velocity gradient tensor.
The rotation tensor is
\begin{equation}
\mathbf{\Omega} = \frac{1}{2} \left( \nabla \mathbf{u} - \nabla \mathbf{u}^\intercal \right).
\end{equation}
The squared norm $||\cdot||^2$ for the tensors is defined as
\begin{equation}
||\mathbf{S}||^2 = \tr{\mathbf{S}^\intercal \mathbf{S}},
\end{equation}
for the symmetric strain tensor and
\begin{equation}
||\mathbf{\Omega}||^2 = \tr{\mathbf{\Omega}^\intercal \mathbf{\Omega}},
\end{equation}
for the anti-symmetric strain tensor.

In regions where \(Q > 0\), the local rotation dominates over strain, and these regions are identified as vortices~\cite{hunt1987}.
The \(Q\) iso-surfaces were visualized at iso values of $Q=0.001$ in \cref{fig:modes}.
These \(Q\)-contours were made at four selected time values $t_1$, $t_2$, $t_3$, $t_4$ such that $t_1 < t_2 < t_3 < t_4$.
The non-dimensional time values were $t_1 = 1075$,  $t_2 = 1076.5$, $t_3 = 1078$ and $t_4 = 1080$.

Physically, the derivative computation of the singular values is the measure of how sensitive each singular value is to perturbations in the original snapshot matrix. 
Thus, we get three directions of derivatives of the velocities as there are three components of the velocity data $u$ , $v$ and $w$. 
With reference to \Cref{eq:snapshot_matrix_representation}, we can see that
\begin{equation}
    \f{\d \sigma_i}{\d \mb{X}} = 
    \begin{bmatrix}
    \f{\d \sigma_i}{\d {u_{1,1}}} & \f{\d \sigma_i}{\d {u_{1,2}}} & \cdots & \f{\d \sigma_i}{\d {u_{1,n}}}\\
    \f{\d \sigma_i}{\d {u_{2,1}}} & \f{\d \sigma_i}{\d {u_{2,2}}} & \cdots & \f{\d \sigma_i}{\d {u_{2,n}}}\\
    \vdots & \vdots & \vdots &\vdots\\
    \f{\d \sigma_i}{\d {u_{m_0,1}}} & \f{\d \sigma_i}{\d {u_{m_0,2}}} & \cdots & \f{\d \sigma_i}{\d {u_{m_0,n}}}\\
    \makebox[0pt][l]{\hspace{-0.6cm}\rule{0.36\linewidth}{0.4pt}}\\
    \f{\d \sigma_i}{\d {v_{1,1}}} & \f{\d \sigma_i}{\d {v_{1,2}}} & \cdots & \f{\d \sigma_i}{\d {v_{1,n}}}\\
    \f{\d \sigma_i}{\d {v_{2,1}}} & \f{\d \sigma_i}{\d {v_{2,2}}} & \cdots & \f{\d \sigma_i}{\d {v_{2,n}}}\\
    \vdots & \vdots & \vdots &\vdots\\
    \f{\d \sigma_i}{\d {v_{m_0,1}}} & \f{\d \sigma_i}{\d {v_{m_0,2}}} & \cdots & \f{\d \sigma_i}{\d {v_{m_0,n}}}\\
    \makebox[0pt][l]{\hspace{-0.6cm}\rule{0.36\linewidth}{0.4pt}}\\
    \f{\d \sigma_i}{\d {w_{1,1}}} & \f{\d \sigma_i}{\d {w_{1,2}}} & \cdots & \f{\d \sigma_i}{\d {w_{1,n}}}\\
    \f{\d \sigma_i}{\d {w_{2,1}}} & \f{\d \sigma_i}{\d {w_{2,2}}} & \cdots & \f{\d \sigma_i}{\d {w_{2,n}}}\\
    \vdots & \vdots & \vdots &\vdots\\
    \f{\d \sigma_i}{\d {w_{m_0,1}}} & \f{\d \sigma_i}{\d {w_{m_0,2}}} & \cdots & \f{\d \sigma_i}{\d {w_{m_0,n}}}\\
    \end{bmatrix}
,
\label{eq:snapshot_matrix_sens}
\end{equation}
where $\d\sigma / \d u$ is the derivative of the singular value with respect to the $u$-velocity.
Similarly this holds for $\d\sigma / \d v$ and $\d\sigma / \d w$.

From this, it is clear that the singular value derivative is directional in the three dimensioned space it is evaluated for through POD.
The snapshot matrix had $149.30352 \times 10^6$ rows equal to the number of states and $75$ columns equal to the number of time steps as aforementioned and the derivatives of the singular value with respect to this snapshot matrix were computed using the RAD formula shown in \Cref{eq:RAD_formula_singular_sens_real}.

\begin{figure}[H]
\centering
\includegraphics[width=1\textwidth]{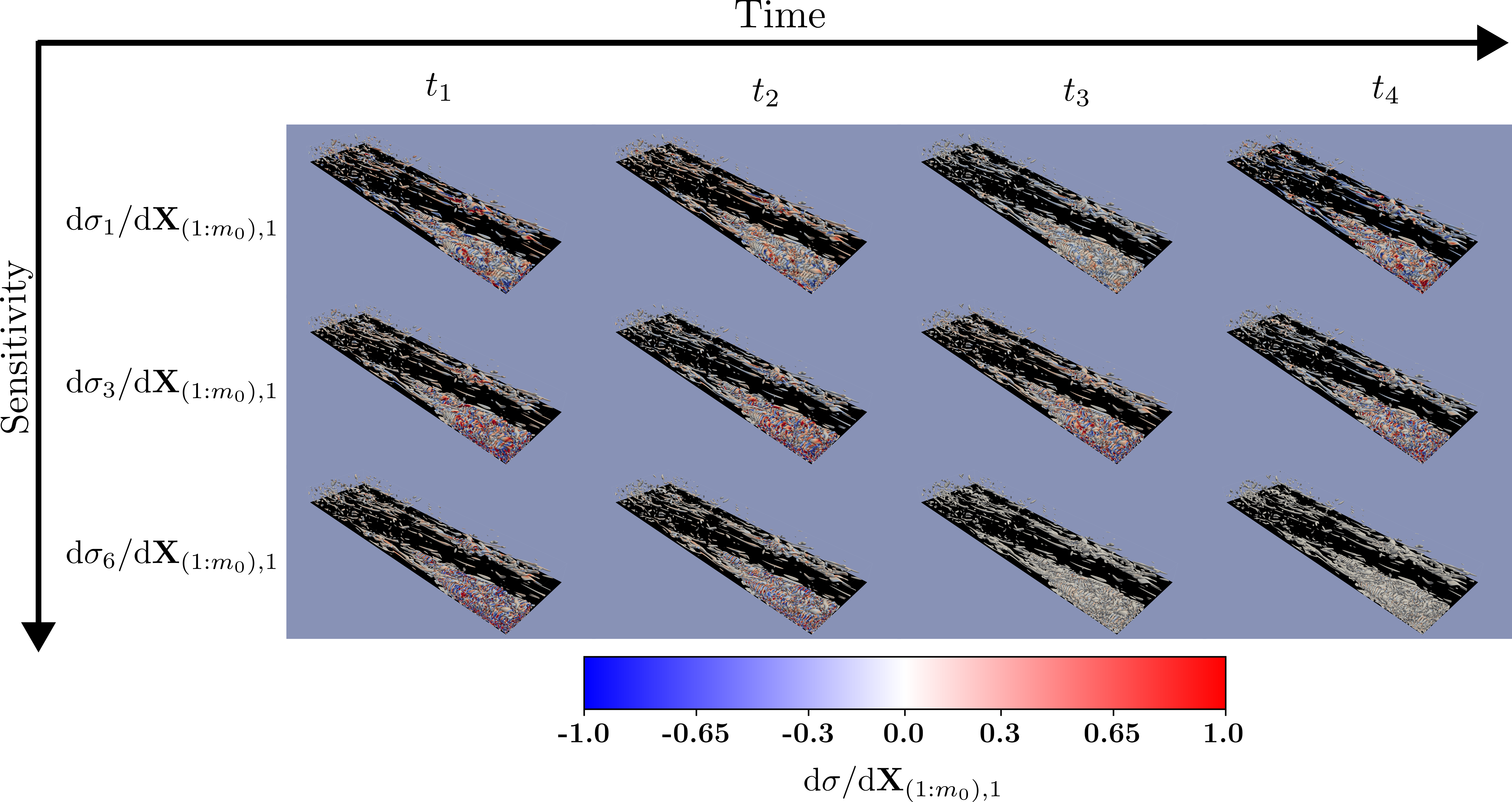}
\caption{$x$-direction singular value derivatives of the POD modes of JHTDB transitional boundary layer dataset. Vortical structures are \(Q\) iso-surfaces (\(Q\) = $0.001$)}
\label{fig:sens}
\end{figure}

\begin{figure}[H]
    \centering
    \includegraphics[width=1\textwidth]{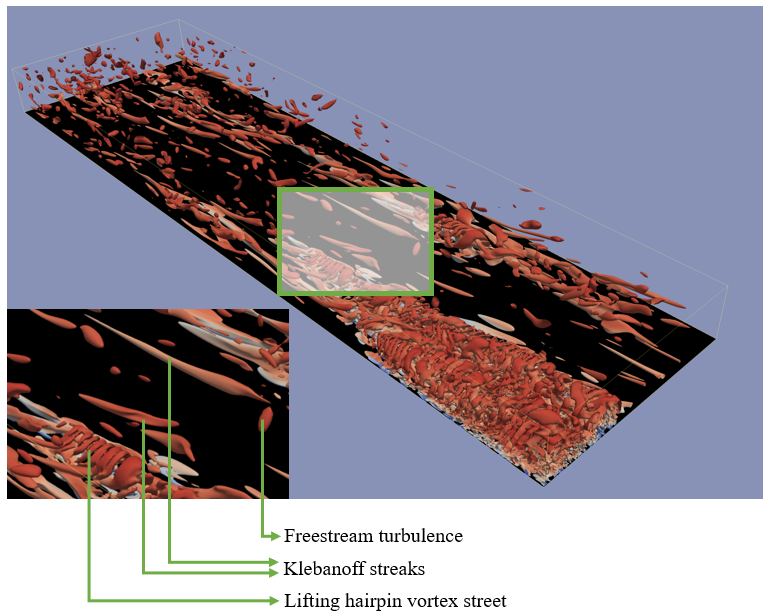}
    \caption{\hlo{Key vortical structures in the transitional region of the flow over flat plate visualized by} (\(Q\) = $0.001$) \hlo{iso-surfaces colored by velocity-magnitude}}
    \label{fig:flow_structs}
\end{figure}

Results from this derivative computation are shown in \cref{fig:sens} which showcases the derivatives of the first, third and sixth modes along the $x$-direction for the first time-step.
In this figure, the $Q$-criterion plot at $Q=0.001$ for the FOM was made and colored by the singular value derivatives.
These $Q$-contours were again made for four selected time values $t_1$, $t_2$, $t_3$, $t_4$ such that $t_1 < t_2 < t_3 < t_4$.
The non-dimensional time values were $t_1 = 1075$,  $t_2 = 1076.5$, $t_3 = 1078$ and $t_4 = 1080$.
\hlo{The following observations were made regarding the flow physics and derivatives of the singular values with respect to the snapshot matrix.
The flow over flat plate case shown is an example of bypass transition.
In typical transition flow, Tollmien--Schlichting waves (TS waves) are seen when the freestream turbulence is negligible~\cite{Schmid2001}.}
\hlo{In the case when the free stream turbulence intensity is about $3.5\%$ the transition occurs with formation of Klebanoff streaks as shown in \cref{fig:flow_structs} in a case called the bypass transition, which essentially ``bypasses'' the formation of TS waves~\cite{Zaki2013}.}

% [x] TODO HS-: We should avoid large png figure like this in the future... Compilation is slow...

\hlo{Downstream of these streaks, the hairpin vortices can be seen which are identified in~\cref{fig:flow_structs} too.
Once the ``lift-up'' effect kicks in, the mixing becomes stronger, more chaotic and the flow transitions to fully turbulent flow~\cite{Adrian2007}. % [x] TODO HS-: "->`` check thoroughly!
The current image snapshots of modes and derivatives of the singular values visualized in~\cref{fig:modes,fig:sens} are shown for the transitional region, which includes the Klebanoff streaks and the hairpin vortices just before the flow transitions to fully turbulent flow~\cite{Lee2018}.}

\hlo{The derivative computation results reveal that the flow structures right before the flow turns fully turbulent are more sensitive to perturbations than the larger and elongated Klebanoff streaks which appear at the start of the transitioning flow.
This suggests that the flow receptivity increases as the flow approaches the fully turbulent regime.
The derivatives fluctuate just as the singular values do with time, which is obvious as the modes are scaled by the singular values along with the temporal coefficients which typically portray a sinusoidal like variation.}

% [x] TODO HS-: incorrect reference to figures. Do a thorigh check.
\hlo{Comparing the identified vortical structures as in} \cref{fig:flow_structs} \hlo{with the structures in}~\cref{fig:sens}\hlo{, it can be seen that the Klebanoff streaks and larger vortical structures at the start of the domain as the flow enters to flow over the flat plate, the structures do not show much receptivity to perturbations in the flow.
However, the smaller structures like the free stream turbulence as well as the hairpin vortices seem to be more sensitive to perturbations.}

\hlo{However, in the temporal sense, it can be seen in}~\cref{fig:sens}\hlo{there comes snapshots in time when there is little to no receptivity in the flow for the smaller modes.
In fact the receptivity of the modes changes periodicially with time as do the temporal coefficients.
The variations are much larger for the smaller modes however, as compared to the more dominant ones.}
% [x] TODO HS-: mention which figure to look at.

\hlo{This kind of identification of sensitive regions of the modes in the flow field offers a pathway for further exploration of receptivity theory and transition to turbulence. 
This work serves as a foundation for future studies aimed at controlling or predicting flow transitions, with potential applications in aerodynamic design, turbulence modeling, and flow control strategies.
This is the general spatial trend observed.}

\hlone{The derivative of singular value of large dataset can be taken using the RAD formula and was shown in this section.
In lieu of this, an in-depth analysis of the flow physics here should be sought, but is deemed as future work by the authors and is therefore beyond the scope of the current manuscript.}

\subsection{\hlo{Wing truss compliance optimization using differentiable SVD}}
% [x] TODO HS-: merge the broad reference into the intro section and directly start with truss problem here.
 
\hlone{In the current section we discuss optimization the bar thickness of a 315-bar wing-truss from the common research model~Brooks et al~\cite{Brooks2018a} by optimizing the singular value of the response matrix to a minimum, essentially reducing the compliance of the structure.
In this case, we start with a stiff wing and make the wing more flexible in simpler terms.}

\hlo{We refer the readers to Ersoy and Mugan~\cite{ersoy2002a}, from where this example's physics was sourced in which the design sensitivity was studied using SVD.
They showed in their paper that the stiffness of a static structure is a measure of the leading singular value of the response matrix of the system.}
% [x] TODO HS-: briefly discuss Gauteux differentials
%RK-HS-: Not required it was unnecessary. Removed it.

\hlo{The net static response of a structure is compliance of that structure and this is inversely proportional to its stiffness.
Thus, reduction of the leading singular value of the response matrix of the system reduces the compliance~\cite{ersoy2002a}.
More details on the mathematical considerations of this statement can be found in}~\cref{sec:truss}.
\hlo{Shahabsafa et al.~\cite{Shahabsafa2018a} used the common research model wing in their design optimization to minimize weight.
In the current study, we optimized the net static response of the same common research model wing truss by optimizing the bar thicknesses in minimization of the leading singular value of the response matrix } (Refer to~\Cref{sec:truss} \hlo{for more details on this matrix) subject to the constraint that the net mass of the truss is constant and bounds placed on the bar areas.}

\begin{figure}[H]
    \centering
    \includegraphics[width=1\textwidth]{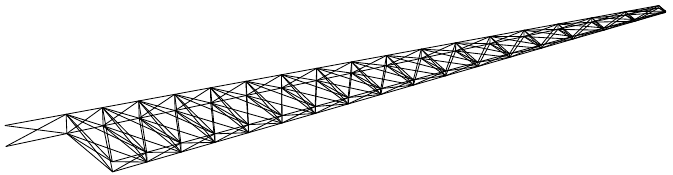}
    \caption{\hlo{The common research model 315 bar wing-truss \cite{Shahabsafa2018a}}}
    \label{fig:wing_truss}
\end{figure}

\hlo{Thus, mass was re-distributed inside the truss amongst the bars and the optimization was carried out.
The motivation of this task is to demostrate the optimization of a wing-truss system using the RAD formula for the derivative computation subject to multiple design variables, which in this case is $315$ since we have $315$ bars.
All values in this optimization such as the material properties and the member aeras were sourced as is from Shahabsafa et al.~\cite{Shahabsafa2018a} and Brooks et al~\cite{Brooks2018a}.}
\hlo{The undeformed wing truss is shown in}~\cref{fig:wing_truss}.

\hlo{We now discuss the derivative computation of the singular value with respect to the design variables -- the bar areas.
The derivative of the first singular value of the response matrix with respect to each of these design variables is a vector given by the matrix-vector product as}
\begin{equation}
    \f{\d \sigma_1}{\d \mb{b}} = \text{vec} \biggl( \f{\d \sigma_1}{\d \mb{K}} \biggl)^\intercal \cdot \biggl[\text{vec} \biggl(\f{\d \mb{K}}{\d b_1} \biggl) \text{~~} \text{vec} \biggl(\f{\d \mb{K}}{\d b_2}\biggl) \cdot\cdot\cdot \text{vec} \biggl(\f{\d \mb{K}}{\d b_n} \biggl) \biggl],
    \label{eq:der_RAD_opt}
\end{equation}
\hlo{where each entry of the vector $\d \sigma_1 / \d \mb{b}$ has the derivative of $\sigma_1$ with respect to the entry's $b_i$, $n$ is the net number of elements of the truss, which is 315 in the undertaken wing truss problem and $\d \sigma_1 / \d \mb{b} \in \mathbb{R}^{n_b}$ where $n_b$ is the number of bars in the truss. %[x] TODO HS-: R^n -> R^n_b? n seems to be used before... 
}
\hlo{The same derivative can also be obtained using FD in which each element of $\mb{b}$ is perturbed to compute the perturbation in $\sigma_1$ and this is computed as $\d \sigma_1 / \d \mb{b}$.}
\hlo{This approach however requires as many derivative function calls as there are design variables.
For each entry of $b_i$, FD requires evaluation of the inverse of stiffness matrix and SVD of this inverse matrix and this becomes computationally intensive for several design variables as was shown in the scalability analysis}~\Cref{sec:scalability}.

\hlo{The optimization problem is formulated as follows}
\begin{equation}
\begin{aligned}
    \min_{\mb{b}} \quad & \sigma_1, \\
    \text{subject to} \quad &  \underline{b} \le b_i \le \bar{b}, \quad i=1, \ldots, n_{\mb{b}},\\
    & \sum_{i=1}^{n} b_i\text{~}l_i = v_0,\\
\end{aligned}
\end{equation}
\hlo{where $l_i$ is the length of the $i$--th bar, $\underline{b}$ is $2.50\times 10^{-5} \text{ m}^2$, $\bar{b}$ is $0.12 \text{ m}^2$, and $v_0=6.79 \text{ m}^3$ is the initial volume. % [x] TODO HS-: check and add units here.
The constraint placed is the volume constraint.
The truss optimization was initialized with the uniform thickness throughout given by area of each element set to $0.01$ $\text{m}^2$.
All material properties, nodes and elemental connections are sourced from Shahabsafa et al.~\cite{Shahabsafa2018a}.}

\begin{figure}[H]
    \centering
    \includegraphics[width=1\textwidth]{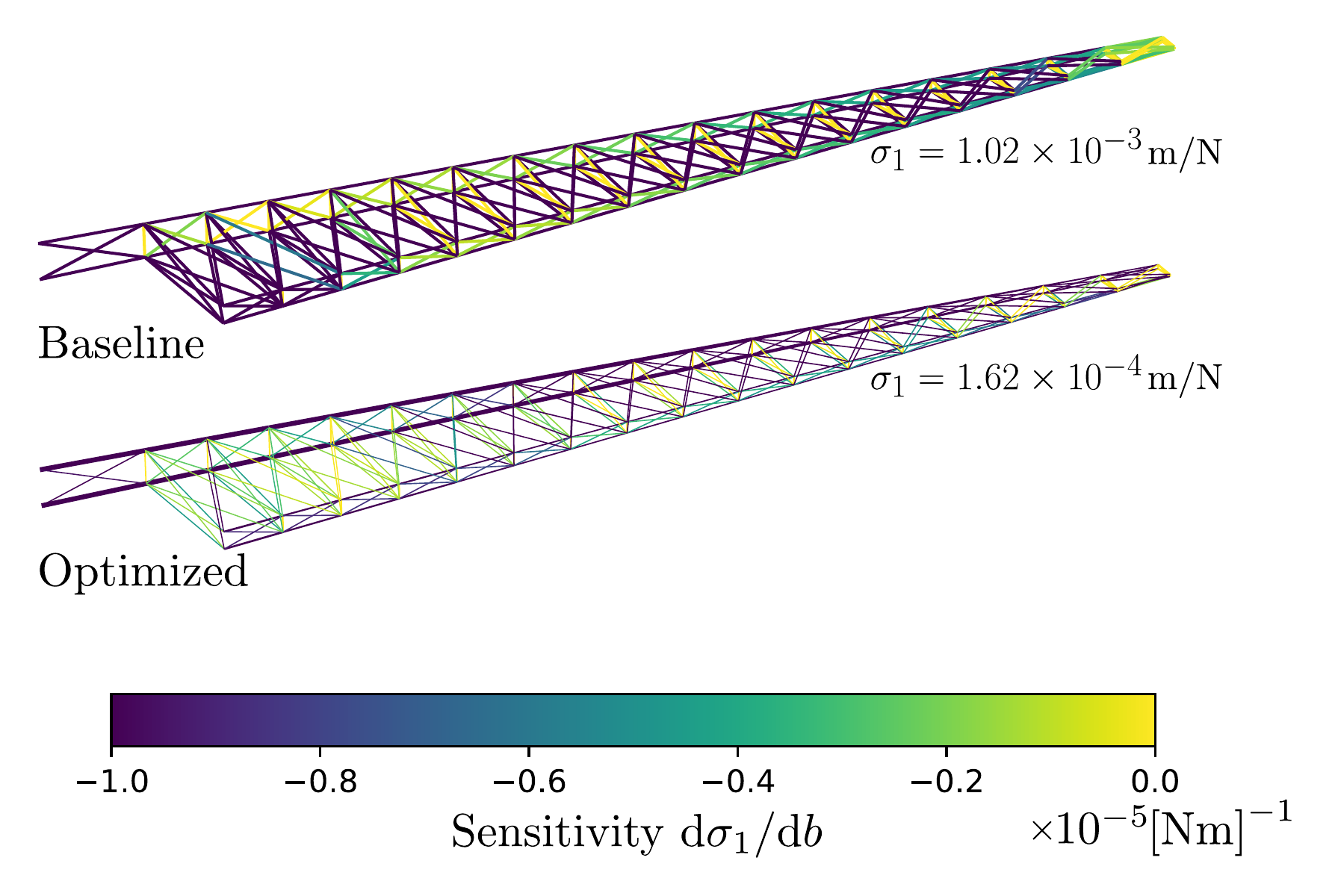}
    \caption{\hlo{The baseline and optimized wing-trusses.
    The bar thickness is represented by the line thickness in the plot.
    Bars are colored by sensitivity magnitude of the singular value.}}
    \label{fig:sens_truss}
\end{figure}
\hlo{The optimization package used was Scipy package's optimizer~\cite{2020SciPy-NMeth} and the optimizer used was based on the trust region methods of gradient based optimization~\cite{Conn2000}.}
% [x] TODO HS-: 1e-5 -> 10^-5 proper latex math format.

% [x] TODO HS-: mentioned optimization package / method
% [x] TODO HS-: add units of variables...
\hlo{The results of optimization are shown in}~\cref{fig:sens_truss}.
\hlo{The objective was minimized with the tolerance of $10^{-6}$ resulting in a final singular value of $1.6 \times 10^{-4} \text{ m} / \text{N}$ from the starting singular value of $1.02 \times 10^{-3} \text{ m} / \text{N}$.
Thus, the net static response of the structure was reduced by $81\%$ post optimization, where the net mass of the truss was re-distributed among the members.}
\hlo{In}~\cref{fig:sens_truss}\hlo{, the thickness of the lines represents the thickness of the bars which was optimized.
Each bar's color represents the leading singular value's sensitivity to the changes in that particular bar's thickness.}

\hlo{Thus, the RAD form of the differentiable SVD approach for real-valued inputs and outputs was used to showcase an optimization of the thickness of the members of a wing truss system to reduce the net static response and the optimization was carried out successfully.}

\section{Conclusions}\label{sec:conc}
In this paper, we developed two adjoint-based methods to compute the singular values and singular vector derivative and an RAD formula to compute the singular value derivative.
The proposed differentiable SVD algorithms do not scale with the number of inputs and do not require all of the singular variables during derivative computation, thus reducing the computational cost, they are relatively easier to implement, and are accurately, being implicit analytic in nature of the solution.
All the derivative computation strategies proposed in the current manuscript can handle complex-valued inputs and outputs.

We proposed two adjoint method approaches leveraging the relationship between EVP and SVD.
In the first approach, we leverage the relationship between the SVD of a matrix and the EVP of its Gram matrices.
We call it the Gram matrix method or GMM.
In the second approach, we leverage the relationship between the SVD of a matrix and the EVP of its symmetric embedding.
We call it the symmetric embedding matrix method or SEMM.
The adjoint-based GMM computes the derivative by breaking down the SVD into two eigenvalue problems of the Gram matrices and applies the adjoint method to the governing equations for each EVP.
This results in two formulations, the Left Gram Matrix Method or LGMM and the Right Gram Matrix Method RGMM, deriving their names from the fact that each Gram matrix either yields the left or right singular vectors in their EVP solutions.

The adjoint-based SEMM computes this derivative by directly applying the adjoint method to the SVD governing equations resulting from the symmetric embedding EVP.
Each method holds its merits based on the objective function whose derivative is desired and the nature of matrix $\mb{A}$ being sparse or dense.
If the matrix is dense, the computation of the Gram matrices of the matrix can become costly.
To avoid this, the SEMM is used.
If the matrix is sparse, the former method, GMM, should be used.
If the matrix is square, LGMM or RGMM are both equivalent.
If the matrix is tall, LGMM is advised to be used, otherwise, RGMM.
Finally, we also proposed a general RAD dot product identity for general complex-valued functions that are imaginary-differentiable.

The results from derivative computations of two imaginary matrices, square and rectangular, were compared against the results from FD approximations, and we achieved a 5--6 digit match.
\hlo{A scalability analysis was carried out to compare and contrast the differences in the computational efficiencies of the proposed methods, RAD formulae in the literature and FD, which resulted in faster computational speeds for the proposed adjoint-RAD methods.}

An implementation of the real-valued form of the proposed singular value derivative RAD formula was shown by applying the formula on the snapshot matrix of a large cutout of the JHTDB dataset in a 3-dimensional space and obtaining the singular values and vectors through POD.
The snapshot matrix had $149.30352 \times 10^6$ rows equal to the number of states and $75$ columns equal to the number of time steps.
\hlo{In another example of wing-truss optimization, the differentiable framework was added as part of an optimization problem of reducing the truss compliance of a 315 bar wing truss source from the undeformed common research model wing.
The optimization was carried out successfully with an $81\%$ reduction in the singular value.}

There is no limitation set on the objective function under consideration, making the proposed method well-suited to large-scale design optimization problems involving gradients of SVD with respect to several design variables, for instance, in differentiable POD or differentiable resolvent analyses that find applications in various engineering-design optimization problems.

\section{Acknowledgement}
The first author would like to thank Ariel Lubonja from John Hopkins University for help with the data retrieval process of the JHTDB flat plate boundary layer dataset.

\bibliographystyle{elsarticle-num-names}
\bibliography{main}

\appendix

\section{Adjoint method for EVP}\label{sec:adjoint_evp}
We derive the adjoint equation for any function, $f = f(\mb{w},\mb{D}(\mb{x}))$, where $\mb{D}$ is the input EVP coefficient matrix and $\mb{w}$ is the solution of the EVP.
The coefficient matrix can be further paramatrized by design variable, $\mb{x}$.
Thus, expanding the adjoint equations \Cref{eq:total_der_eqn} and \Cref{eq:adjoint_eqn} from \cref{sec:Governing_equns_and_adjoint_method}, we have
\begin{equation}
\begin{bmatrix}
\mb{D}_r - \lambda_r \mb{I} & -\mb{D}_i + \lambda_i \mb{I} & -\boldsymbol{\phi}_r & \boldsymbol{\phi}_i\\
\mb{D}_i - \lambda_i \mb{I} & \mb{D}_r + \lambda_r \mb{I} & -\boldsymbol{\phi}_i & -\boldsymbol{\phi}_r\\
2\boldsymbol{\phi}_r^\intercal & 2\boldsymbol{\phi}_i^\intercal & 0 & 0\\
0 & \boldsymbol{e}_k^\intercal & 0 & 0
\end{bmatrix}
\begin{bmatrix}
\boldsymbol{\psi}_{\text{main,r}} \\
\boldsymbol{\psi}_{\text{main,i}} \\
\boldsymbol{\psi}_{m}\\
\boldsymbol{\psi}_{p}
\end{bmatrix}
= \f{\partial f}{\partial \mb{w}}^\intercal,
\label{eq:adjoint_equn_EVP_appendix}
\end{equation}
which is the adjoint equation for the EVP of matrix $\mb{D}$.
The partial derivative of the residual matrix with respect to the imaginary parts of the matrix $\mb{D}$ is
\begin{equation}
\begin{aligned}
\f{\partial \mb{r}}{\partial \mb{D}_r} &= \boldsymbol{\psi}_{\text{main,r}} \boldsymbol{\phi}_r^\intercal + \boldsymbol{\psi}_{\text{main,i}} \boldsymbol{\phi}_i^\intercal,\\
\f{\partial \mb{r}}{\partial \mb{D}_i} &= -\boldsymbol{\psi}_{\text{main,r}} \boldsymbol{\phi}_i^\intercal + \boldsymbol{\psi}_{\text{main,i}} \boldsymbol{\phi}_r^\intercal,
\label{eq:a}
\end{aligned}
\end{equation}
and the partial derivatives of imaginary components of $f$ with respect to those of the matrix $\mb{D}$ are 
\begin{equation}
\begin{aligned}
\f{\d f_r}{\d \mb{D}_r} = - \f{\partial \mb{r}}{\partial \mb{D}_r}^\intercal \boldsymbol{\psi}_r,\\
\f{\d f_r}{\d \mb{D}_i} = - \f{\partial \mb{r}}{\partial \mb{D}_i}^\intercal \boldsymbol{\psi}_r,\\
\f{\d f_i}{\d \mb{D}_r} = - \f{\partial \mb{r}}{\partial \mb{D}_r}^\intercal \boldsymbol{\psi}_i,\\
\f{\d f_i}{\d \mb{D}_i} = - \f{\partial \mb{r}}{\partial \mb{D}_i}^\intercal \boldsymbol{\psi}_i.
\end{aligned}
\end{equation}
This completes the adjoint method for solving EVP. 
Detailed approach can be found in He et al~\cite{He2023}.

\section{Chain rule RAD formula}\label{sec:Chain_rule_RAD_formula}

Consider the matrix $\mb{B}=\mb{A}\mb{A}^*$. An FAD application on both sides results in the real part as
\begin{equation}
\dot{\mb{B}}_r = \dot{\mb{A}}_r \mb{A}_r^\intercal + \mb{A}_r \dot{\mb{A}}_r^\intercal + \dot{\mb{A}}_i\mb{A}_i^\intercal + \mb{A}_i \dot{\mb{A}}_i^\intercal,
\end{equation}
and the imaginary part as
\begin{equation}
\dot{\mb{B}}_i = \dot{\mb{A}}_i \mb{A}_r^\intercal + \mb{A}_i \dot{\mb{A}}_r^\intercal - \dot{\mb{A}}_r\mb{A}_i^\intercal - \mb{A}_r \dot{\mb{A}}_i^\intercal.
\end{equation}
These two terms $\dot{\mb{B}}_r$ and $\dot{\mb{B}}_i$ can be used in the Trace expression in \Cref{eq:trace_identity_proposed} as
\begin{equation}
\mathrm{Tr}[\overbar{\mb{B}}_r^\intercal\dot{\mb{B}}_r + \overbar{\mb{B}}_i^\intercal\dot{\mb{B}}_i] = \mathrm{Tr}[\overbar{\mb{A}}_r^\intercal\dot{\mb{A}}_r + \overbar{\mb{A}}_i^\intercal\dot{\mb{A}}_i].
\label{eq:Tr_identity_applied_to_AA*}
\end{equation}
Upon applying the Tr identities and expanding LHS of \Cref{eq:Tr_identity_applied_to_AA*}, we have
\begin{equation}
\begin{aligned}
\mathrm{Tr}[(\mb{A}_r^\intercal \overbar{\mb{B}}_r^\intercal + \mb{A}_r^\intercal \overbar{\mb{B}}_r + \mb{A}_i^\intercal \overbar{\mb{B}}_i - \mb{A}_i^\intercal \overbar{\mb{B}}_i^\intercal)\dot{\mb{A}}_r + \\
(\mb{A}_i^\intercal \overbar{\mb{B}}_r^\intercal + \mb{A}_i^\intercal \overbar{\mb{B}}_r + \mb{A}_r^\intercal \overbar{\mb{B}}_i^\intercal - \mb{A}_r^\intercal \overbar{\mb{B}}_i)\dot{\mb{A}}_i] \\
= \mathrm{Tr}[\overbar{\mb{A}}_r^\intercal\dot{\mb{A}}_r + \overbar{\mb{A}}_i^\intercal\dot{\mb{A}}_i],
\end{aligned}
\label{eq:Tr_identity_applied_to_AA*_expanded}
\end{equation}
which upon simplifying we get for the real part
\begin{equation}
\overbar{\mb{A}}_r = {(\mb{A}_r^\intercal \overbar{\mb{B}}_r^\intercal + \mb{A}_r^\intercal \overbar{\mb{B}}_r + \mb{A}_i^\intercal \overbar{\mb{B}}_i - \mb{A}_i^\intercal \overbar{\mb{B}}_i^\intercal)}^\intercal,
\label{eq:Ar_bar_AA*}
\end{equation}
and for the imaginary part
\begin{equation}
\overbar{\mb{A}}_i = {(\mb{A}_i^\intercal \overbar{\mb{B}}_r^\intercal + \mb{A}_i^\intercal \overbar{\mb{B}}_r + \mb{A}_r^\intercal \overbar{\mb{B}}_i^\intercal - \mb{A}_r^\intercal \overbar{\mb{B}}_i)}^\intercal.
\label{eq:Ai_bar_AA*}
\end{equation}
Now we can place the reverse seed as $g_r$ and $g_i$ in each of the \Cref{eq:Ar_bar_AA*} and \Cref{eq:Ai_bar_AA*} to obtain the derivatives as in \Cref{eq:get_dfdA_from B}.

Consider then the matrix $\mb{C}=\mb{A}^*\mb{A}$
Following the same methods as in \Cref{eq:Tr_identity_applied_to_AA*,eq:Tr_identity_applied_to_AA*_expanded,eq:Ar_bar_AA*,eq:Ai_bar_AA*}, we have for the real part
\begin{equation}
\overbar{\mb{A}}_r = {(\overbar{\mb{C}}_r \mb{A}_r^\intercal + \overbar{\mb{C}}_r^\intercal \mb{A}_r^\intercal + \overbar{\mb{C}}_i \mb{A}_i^\intercal - \overbar{\mb{C}}_i^\intercal \mb{A}_i^\intercal)}^\intercal ,
\label{eq:Ar_bar_A*A}
\end{equation}
and for the imaginary part
\begin{equation}
\overbar{\mb{A}}_i = {(\overbar{\mb{C}}_r \mb{A}_i^\intercal + \overbar{\mb{C}}_r^\intercal \mb{A}_i^\intercal + \overbar{\mb{C}}_i \mb{A}_r^\intercal - \overbar{\mb{C}}_i^\intercal \mb{A}_r^\intercal)}^\intercal .
\label{eq:Ai_bar_A*A}
\end{equation}
We can then seed $h_r$ and $h_i$ each into \Cref{eq:Ar_bar_A*A,eq:Ai_bar_A*A} to obtain the derivatives in \Cref{eq:get_dfdA_from_C}.

\section{Vectorization}
\label{sec:vectorization}

The vectorization operator $\vop\left(\cdot\right)$ is defined as
\begin{equation}
\left(\vop\left(\mb{A}\right)\right)_{i\times (n_2-1) + j} = \mb{A}_{ij}, \quad i=1, \ldots, n_1, \quad j=1, \ldots, n_2,
\end{equation}
where $\mb{A} \in \mathbb{R}^{n_1\times n_2}$, $\vop: \mathbb{R}^{n_1\times n_2} \rightarrow \mathbb{R}^{n_1n_2}$, and the subscript represents the index of an element from the matrix $\mb{A}$ or the vector $\vop\left(\mb{A}\right)$.
This linear operation transforms a matrix into a vector simplifying the matrix derivative computation.
The inverse vectorization operator $\vop^{-1}\left(\cdot\right)$ is defined as,
\begin{equation}
\vop^{-1}\left(\vop\left(\mb{A}\right)\right)= \mb{A},
\end{equation}
for arbitrary matrix $\mb{A}$.
As a special case of the properties of vectorization,
\begin{equation}
\begin{aligned}
\vop\left(\mb{a}\right) &= \mb{a},\\
\vop^{-1}\left(\mb{a}\right) &= \mb{a},
\end{aligned}
\end{equation}
where $\mb{a} \in \mathbb{R}^{n_3}$ is some arbitrary vector.

The following convention is used when writing a derivative involving matrices in this paper for  $\left(\p \mb{A} / \p \mb{B}\right)^\intercal \overbar{\mb{A}}$, where $\mb{A}\in \mathbb{R}^{n_1\times n_2},$ and $\mb{B}\in\mathbb{R}^{n_2\times n_3}$.
Using the vectorization notation, $\left(\p \mb{A} / \p \mb{B}\right)^\intercal \overbar{\mb{A}}$ is a simplified notation of
\begin{equation}
\vop^{-1}\left(\left(\f{\p \vop\left(\mb{A}\right)}{\p \vop\left(\mb{B}\right)}\right)^\intercal \vop\left(\overbar{\mb{A}}\right)\right) \in \mathbb{R}^{n_2\times n_3}.
\end{equation}
This signifies the computation of the tensor-vector products such as seen in \Cref{eq:get_dfdA_from B,eq:get_dfdA_from_C}.

\section{RAD form for GMM}\label{sec:RAD_form_for_GMM}

Consider the governing equation for SVD as shown in \Cref{eq:SVD_gov_eqn}.
The vector $\mb{u}$ can be written as 
\begin{equation}
\mb{u} = \f{\mb{Av}}{\sigma},
\end{equation}
which can be re-written in terms of the real and imaginary parts as 
\begin{equation}
\begin{aligned}
\mb{u}_r = \mb{A}_r \mb{v}_r \f{1}{\sigma} - \mb{A}_i \mb{v}_i \f{1}{\sigma},\\
\mb{u}_i = \mb{A}_r \mb{v}_i \f{1}{\sigma} - \mb{A}_i \mb{v}_r \f{1}{\sigma}.\\
\end{aligned}
\end{equation}
Applying the dot product identity from~\Cref{sec:dot-product-identities} and considering the contribution from $\mb{A}$, we have
\begin{equation}
\begin{aligned}
\overbar{\mb{A}}_r = \f{1}{\sigma}\bigg(\overbar{\mb{u}}_r \mb{v}_r^\intercal + \overbar{\mb{u}}_i \mb{v}_i^\intercal \bigg),\\
\overbar{\mb{A}}_i = \f{1}{\sigma}\bigg(\overbar{\mb{u}}_i \mb{v}_r^\intercal - \overbar{\mb{u}}_r \mb{v}_i^\intercal \bigg).\\
\end{aligned}
\end{equation}
Consider then the same governing equation for SVD as described in \Cref{eq:SVD_gov_eqn}.
For the second route taken for the Adjoint-based derivative computation described in \cref{sec:GMM_sens}, the vector $\mb{v}$ is expressed in terms of the remaining singular variables as 
\begin{equation}
\mb{v} = \f{\mb{A}^*\mb{u}}{\sigma}.
\end{equation}
Thus, through trace identities from~\Cref{sec:dot-product-identities}, we have 
\begin{equation}
\begin{aligned}
\overbar{\mb{A}}_r = \f{1}{\sigma}\bigg(\mb{u}_r \overbar{\mb{v}}_r^\intercal + \mb{u}_i \overbar{\mb{v}}_i^\intercal \bigg),\\
\overbar{\mb{A}}_i = \f{1}{\sigma}\bigg(\mb{u}_i \overbar{\mb{v}}_r^\intercal - \mb{u}_r \overbar{\mb{v}}_i^\intercal \bigg).\\
\end{aligned}
\end{equation}
This concludes the derivation of for the term in the brackets of \Cref{eq:get_dfdA_from B,eq:get_dfdA_from_C}.

\section{Derivation of \Cref{eq:get_dfdA_components_direct_SVD}}\label{sec:derivation}

Consider the residual form $\mb{r}(\mb{w})$ in \Cref{eq:residual_form_svd_gov_eqn}.
Upon differentiating this with respect to $\mb{A}_r$, we have in FAD notation
\begin{equation}
\dot{\mb{r}} = 
\begin{bmatrix}
\dot{\mb{A}_r} \mb{v}_r \\
\dot{\mb{A}_r} \mb{v}_i \\
\dot{\mb{A}_r}^\intercal \mb{u}_r \\
\dot{\mb{A}_r}^\intercal \mb{u}_i \\
0 \\
0
\end{bmatrix}.
\label{eq:dr_wrt_Ar}
\end{equation}
Applying the dot product identity for real valued functions here, we have
\begin{equation}
\mathrm{Tr}(\overbar{\mb{r}}^\intercal \dot{\mb{r}}) = \mathrm{Tr}(\overbar{\mb{A}}_r^\intercal \dot{\mb{A}}_r).
\end{equation}
Placing the adjoint vector $\boldsymbol{\psi}$ as a reverse seed here, we get
\begin{equation}
\mathrm{Tr}(\overbar{\mb{r}}^\intercal \dot{\mb{r}}) = \mathrm{Tr}(\boldsymbol{\psi}^\intercal \dot{\mb{r}}),
\end{equation}
which upon expanding the RHS, we get the expression
\begin{equation}
\mathrm{Tr}(\boldsymbol{\psi}_{v_r}^\intercal \dot{\mb{A}}_r \mb{v}_r + \boldsymbol{\psi}_{v_i}^\intercal \dot{\mb{A}}_r \mb{v}_i + \boldsymbol{\psi}_{u_r}^\intercal \dot{\mb{A}}_r^\intercal \mb{u}_r + \boldsymbol{\psi}_{u_i}^\intercal \dot{\mb{A}}_r^\intercal \mb{u}_i),
\end{equation}
and then applying the trace identities and re-arranging terms, we get
\begin{equation}
\mathrm{Tr}( (\mb{v}_r \boldsymbol{\psi}_{v_r} ^\intercal + \mb{v}_i \boldsymbol{\psi}_{v_i} ^\intercal + \boldsymbol{\psi}_{u_r} \mb{u}_r ^\intercal + \boldsymbol{\psi}_{u_i} \mb{u}_i ^\intercal ) \dot{\mb{A}}_r) = \mathrm{Tr}(\overbar{\mb{A}}_r^\intercal \dot{\mb{A}}_r).
\end{equation}
Since this must hold for any arbitrary $\dot{\mb{A}}_r$, we have
\begin{equation}
\overbar{\mb{A}}_r = \boldsymbol{\psi}_{v_r} \mb{v}_r ^\intercal + \boldsymbol{\psi}_{v_i} \mb{v}_i ^\intercal + \mb{u}_r \boldsymbol{\psi}_{u_r} ^\intercal + \mb{u}_i \boldsymbol{\psi}_{u_i} ^\intercal.
\end{equation}
Similarly taking the derivative of $\mb{r}(\mb{w})$ in \Cref{eq:residual_form_svd_gov_eqn} with respect to $\mb{A}_i$ and following a similar procedure as described above, we have
\begin{equation}
\overbar{\mb{A}}_i = -\boldsymbol{\psi}_{v_r} \mb{v}_i ^\intercal + \boldsymbol{\psi}_{v_i} \mb{v}_r ^\intercal + \mb{u}_i \boldsymbol{\psi}_{u_r} ^\intercal - \mb{u}_r \boldsymbol{\psi}_{u_i} ^\intercal.
\end{equation}
This concludes our derivation of \Cref{eq:get_dfdA_components_direct_SVD}.
It must be noted that the derivative here is  ${(\p \mb{r} / \p \mb{A})}^\intercal \boldsymbol{\psi}$.
It is obtained by reverse seeding the residual vector $\mb{r}$.

\section{Singular value derivative using RAD}\label{sec:Singular_Val_Sens_RAD_Derivation}
Consider the singular value (dominant) written in terms of the imaginary matrix, $\mb{A}$, and its corresponding singular vectors as
\begin{equation}
\label{eq:singular_value_expression_appendix}
\sigma = \mb{u}^* \mb{A} \mb{v},
\end{equation}
where ``$\square^*$'' represents the imaginary conjugate operation. 
$\mb{u}$ and $\mb{v}$ are the left and right singular vectors respectively for the singular value.
We will now derive the derivative of this singular value with respect to the matrix $\mb{A}$. 
\Cref{eq:singular_value_expression_appendix} can be re-written as 

\begin{equation}
\begin{aligned}
\mb{A}\mb{v} &= \sigma \mb{u},\\
\Rightarrow \dot{\mb{A}} \mb{v} + \mb{A}\dot{\mb{v}} &=\dot{\sigma} \mb{u} + \sigma \dot{\mb{u}},\\
\Rightarrow \mb{u}^*\dot{\mb{A}} \mb{v} + \mb{u}^*\mb{A}\dot{\mb{v}} &=\mb{u}^*\dot{\sigma} \mb{u} + \sigma \mb{u}^*\dot{\mb{u}},\\
\Rightarrow \mb{u}^*\dot{\mb{A}} \mb{v} + \mb{u}^*\mb{U}\mb{\Sigma}\mb{V}^*\dot{\mb{v}} &=\mb{u}^*\dot{\sigma} \mb{u} + \sigma \mb{u}^*\dot{\mb{u}},\\
\Rightarrow \mb{u}^*\dot{\mb{A}} \mb{v} + \sigma\mb{v}^*\dot{\mb{v}} &=\mb{u}^*\dot{\sigma} \mb{u} + \sigma \mb{u}^*\dot{\mb{u}},\\
\Rightarrow \mb{u}^*\dot{\mb{A}} \mb{v} &=\dot{\sigma},\\
\end{aligned}
\label{eq:FAD_Form_singular_val_RAD}
\end{equation}
where going from the second from last to the last equation, we used the fact that $\mb{u}^*\mb{u}=\mb{v}^*\mb{v}=1$, thus, $\dot{\mb{u}}^*\mb{u}=\dot{\mb{v}}^*\mb{v}=0$.

Expanding the final equation in \Cref{eq:FAD_Form_singular_val_RAD} into its real and imaginary parts, we have
\begin{equation}
\begin{aligned}
\dot{\sigma} =\text{~} &\mb{u}_r^\intercal \dot{\mb{A}}_r \mb{v}_r + \mb{u}_i^\intercal \dot{\mb{A}}_i \mb{v}_r - \mb{u}_r^\intercal \dot{\mb{A}}_i \mb{v}_i + \mb{u}_i^\intercal \dot{\mb{A}}_r \mb{v}_i \\
&+ i( \mb{u}_r^\intercal \dot{\mb{A}}_i \mb{v}_r - \mb{u}_i^\intercal \dot{\mb{A}}_r \mb{v}_r + \mb{u}_r^\intercal \dot{\mb{A}}_r \mb{v}_i + \mb{u}_i^\intercal \mb{A}_i \mb{v}_i),
\end{aligned}
\label{eq:FAD_Form_singular_val_RAD_real_imag}
\end{equation}
where the imaginary part is zero because the singular value is a real number.
Thus,
\begin{equation}
\begin{aligned}
\dot{\sigma} =\text{~} &\mb{u}_r^\intercal \dot{\mb{A}}_r \mb{v}_r + \mb{u}_i^\intercal \dot{\mb{A}}_i \mb{v}_r - \mb{u}_r^\intercal \dot{\mb{A}}_i \mb{v}_i + \mb{u}_i^\intercal \dot{\mb{A}}_r \mb{v}_i,
\end{aligned}
\label{eq:FAD_Form_singular_val_RAD_real}
\end{equation}
which is the final FAD form to begin with for applying the dot product identity formula from Giles~\cite{Giles2008} as
\begin{equation}
\begin{aligned}
\tr{\overbar{\sigma}^\intercal \dot{\sigma}} = \tr{\overbar{\sigma}^\intercal (\mb{u}_r^\intercal \dot{\mb{A}}_r \mb{v}_r + \mb{u}_i^\intercal \dot{\mb{A}}_i \mb{v}_r - \mb{u}_r^\intercal \dot{\mb{A}}_i \mb{v}_i + \mb{u}_i^\intercal \dot{\mb{A}}_r \mb{v}_i)},
\end{aligned}
\end{equation}
Applying trace identities, we can re write it for $\mb{A}_r$ as
\begin{equation}
\tr{\overbar{\sigma}^\intercal (\mb{v}_r \mb{u}_r^\intercal + \mb{v}_i \mb{u}_i^\intercal)\dot{\mb{A}_r}} = \tr{\overbar{\mb{A}_r}^\intercal \dot{\mb{A}_r}}.
\end{equation}
This has to hold for any arbitrary $\dot{\mb{A}_r}$. 
Thus, we have
\begin{equation}
\overbar{\mb{A}}_r = (\mb{u}_r \mb{v}_r^\intercal + \mb{u}_i \mb{v}_i^\intercal)\overbar{\sigma}.
\label{eq:Ar_bar_sigma_appendix}
\end{equation}
Similarly for $\mb{A}_i$, we have:
\begin{equation}
\tr{\overbar{\sigma}^\intercal (\mb{v}_r \mb{u}_i^\intercal - \mb{v}_i \mb{u}_r^\intercal)\dot{\mb{A}_r}} = \tr{\overbar{\mb{A}_i}^\intercal \dot{\mb{A}_i}}.
\end{equation}
This has to hold for any arbitrary $\dot{\mb{A}_i}$. 
Thus, we have
\begin{equation}
\overbar{\mb{A}}_i = (-\mb{u}_r \mb{v}_i^\intercal + \mb{u}_i \mb{v}_r^\intercal)\overbar{\sigma}.
\label{eq:Ai_bar_sigma_appendix}
\end{equation}
Placing the reverse seed of $\overbar{\sigma}$ equal to 1 in \Cref{eq:Ar_bar_sigma_appendix} and \Cref{eq:Ai_bar_sigma_appendix}, we can obtain the partial derivatives $\d \sigma_r/\d \mb{A}_r$ and $\d \sigma_r/\d \mb{A}_i$. 

In the real dimension, we will have upon placing the reverse seed of $\overbar{\sigma}$ equal to 1 in \Cref{eq:Ar_bar_sigma_appendix},
\begin{equation}
\f{\d \sigma}{\d \mb{A}} = \mb{u} \mb{v}^\intercal,
\label{eq:appendix_RAD_sigma}
\end{equation}
where $\mb{u} \in \mathbb{R} $ and $\mb{v} \in \mathbb{R}$.
This concludes the derivation for \Cref{eq:RAD_formula_singular_sens_imaginary}.

\section{Adjoint method Jacobian computation and direct contribution of main matrix}\label{sec:Jacobian_and_JAX}
Consider the objective function $f$ as 
\begin{equation}
f = \mb{c}^\intercal \mb{u},
\label{eq:ex_func}
\end{equation}
where $\mb{c}$ is a imaginary constant vector and $\mb{u}$ is the left singular vector.
The function $f$ can be written in its real and imaginary parts as $f = f_r + i f_i$.
The real and imaginary parts are  
\begin{equation}
\begin{aligned}
f_r = \mb{c}_r^\intercal \mb{u}_r - \mb{c}_i^\intercal \mb{u}_i ,\\
f_i = \mb{c}_r^\intercal \mb{u}_i + \mb{c}_i^\intercal \mb{u}_r.
\end{aligned}
\end{equation}
Here, the Jacobian of the components of $f$ ($f_r$ and $f_i$) with respect to the state variables $\mb{w}$ for the SEMM case is  $\p f / \p \mb{w}$ where $\mb{w}$ can be taken from \Cref{eq:residual_form_svd_gov_eqn}.
The Jacobian for $f_r$ will then look like 
\begin{equation}
\f{\p f_r}{\p \mb{w}} = 
\begin{bmatrix}
\mb{c}_r\\
-\mb{c}_i\\
0\\
0\\
0\\
0\\
\end{bmatrix}.
\end{equation}
Similarly the Jacobian for $f_i$ will look like 
\begin{equation}
\f{\p f_i}{\p \mb{w}} = 
\begin{bmatrix}
\mb{c}_i\\
\mb{c}_r\\
0\\
0\\
0\\
0\\
\end{bmatrix}.
\end{equation}
Once we have these Jacobians, we can proceed with the derivative computation method discussed in \cref{sec:Numerical_results}, particularly from \Cref{eq:SVD_adjoint_eqn} since we now have its RHS.
A similar approach can be shown for the \Cref{eq:Adjoint_eqn_for_g,eq:Adjoint_eqn_for_h} in the GMM case for formation of the Jacobian and hence the RHS of the adjoint equations.
It can be seen however that this is simply cumbersome even for such a simple linear function as shown in \Cref{eq:ex_func}.

Thus, JAX~\cite{jax2018github} was used for this step, and so can any other equivalent automatic differentiation tool be used to save implementation effort and time.
The same approach can be taken for the direct dependencies of the function on the matrix $\mb{A}$, if any.
If the dependency arises from the GMM computation, then the formulae shown in~\Cref{sec:RAD_form_for_GMM} can be used, as pointed out in \cref{sec:GMM_sens}.
These formulae are based on the automatic differentiation and hence are equivalent when either this approach or the JAX approach is used for evaluation of this step.

\section{FD formula}\label{sec:FD_formula}
A finite difference formula (FD) was employed to compare the results from the proposed Adjoint and RAD methods for the singular value derivative in \cref{sec:Square_matrix_case} and \cref{sec:Rect_mat_case}.
The following equations were employed for the FD results in \cref{tab:Adjoint_FD_1,tab:Adjoint_FD_2,tab:RAD_sing_sens_table_Square,tab:RAD_sing_sens_table_rect}.
\begin{equation}
\begin{aligned}
\f{\d f_r}{\d \mb{A}_{\text{r,pq}}} &= \text{Re}\biggl(\f{f(\mb{A} + \epsilon \mb{E}_{\text{pq}}) - f(\mb{A})}{\epsilon} \biggl),\\
\f{\d f_r}{\d \mb{A}_{\text{i,pq}}} &= \text{Im}\biggl(\f{f(\mb{A} + i\epsilon \mb{E}_{\text{pq}}) - f(\mb{A})}{\epsilon} \biggl),\\
\f{\d f_i}{\d \mb{A}_{\text{r,pq}}} &= \text{Re}\biggl(\f{f(\mb{A} + \epsilon \mb{E}_{\text{pq}}) - f(\mb{A})}{\epsilon} \biggl),\\
\f{\d f_i}{\d \mb{A}_{\text{i,pq}}} &= \text{Im}\biggl(\f{f(\mb{A} + i\epsilon \mb{E}_{\text{pq}}) - f(\mb{A})}{\epsilon} \biggl),
\end{aligned}
\label{eq:FD_equns}
\end{equation}
where $\epsilon = 10^{-6}$ is the finite perturbation applied to the matrix $\mb{A}$, $\mb{E}_{\text{pq}}$ is a single entry matrix with the element on $\text{p}^\text{th}$ row and $\text{q}^\text{th}$ column set equal to $\epsilon$ and the remaining elements set to zero.

\section{The method of snapshots}\label{sec:mthod_of_snapshots}

In the method of snapshots, a covariance matrix is formed using the equation
\begin{equation}
\mb{C} = \mb{X}^\intercal \mb{X},
\end{equation}
where $\mb{X}$ is the snapshot data matrix.
The EVP of $\mb{X}$ is then solved to get the eigenvectors and eigenvalues
\begin{equation}
\mb{C}\mb{v}_i = \lambda_i \mb{v}_i.
\end{equation}
The POD modes are then simply
\begin{equation}
\mb{\Phi}_i = \mb{X} \mb{v}_i \f{1}{\sqrt{\lambda_i}},
\end{equation}
where $i$ denotes each mode level.
The energy levels are related to $\lambda$, and the singular values from the SVD are simply the square roots of all the $\lambda_i$.

The eigenvectors $\mb{v}_i$ are the same as the right singular vectors of SVD of $\mb{X}$ and are related to the temporal coefficients.
The temporal coefficients can be calculated as
\begin{equation}
a_i(t) = \mb{v}_i\sqrt{\lambda_i},
\end{equation}
where $t$ is the physical time step value.

\section{RAD dot product identities and chain rule}\label{sec:dot-product-identities}
We propose the dot product identity for imaginary differentiable function (in multivariate sense).
For real valued matrices, $\mb{A}\in\mathbb{R}^{m_1\times n_1}$, and functions of those matrices, $\mb{B}(\mb{A}):\mathbb{R}^{m_1\times n_1} \rightarrow \mathbb{R}^{m_2\times n_2}$, \citet{Giles2008} proposed the following dot product identity
\begin{equation}
\tr{\overbar{\mb{B}}^\intercal \dot{\mb{B}}} = \tr{(\overbar{\mb{A}}^\intercal \dot{\mb{A}}},
\label{eq:Tr_identity_real}
\end{equation}
where $\mb{B}$ is a real differentiable function of $\mb{A}$. 

He et al~\cite{He2023} proposed the following identity for imaginary matrices, $\mb{A}\in\mathbb{C}^{m_1\times n_1}$, and analytic matrix functions of those matrices, $\mb{B}(\mb{A}):\mathbb{C}^{m_1\times n_1}\rightarrow \mathbb{C}^{m_2\times n_2}$
\begin{equation}
\tr{\overbar{\mb{B}}^* \dot{\mb{B}}} = \tr{\overbar{\mb{A}}^* \dot{\mb{A}}},
\label{eq:trace_identity_analytic}
\end{equation}
where the tranpose operator, $\square^\intercal$, is replaced by the conjugate transpose operator, $\square^*$.

However, for imaginary differentiable function (in multivariate sense) but not imaginary non-analytic functions, the conjugate variables need to be considered. 
After applying the Wiritinger derivatives, Roberts and Roberts proposed the following RAD trace-identity for imaginary matrices and their non-analytic functions:
\begin{equation}
\label{eq:trace_identity_non_analytic}
\tr{\overbar{\mb{B}}^* \dot{\mb{B}} + {(\overbar{\mb{B}}^* \dot{\mb{B}})}^*} = \tr{\overbar{\mb{A}}^* \dot{\mb{A}} + {(\overbar{\mb{A}}^* \dot{\mb{A}})}^*}.
\end{equation}
Notice that for $\mb{B}$ being an analytic function of $\mb{A}$, the second term in the Trace operator on each side of \Cref{eq:trace_identity_non_analytic} becomes zero, effectively reducing to \Cref{eq:trace_identity_analytic}.

We can also expand \Cref{eq:trace_identity_non_analytic} into real and imaginary components which is more convenient to deal with in many cases where the real and imaginary components are derived separately.
\Cref{eq:trace_identity_non_analytic} can be written as
\begin{equation}
\label{eq:trace_identity_proposed}
\mathrm{Tr}\left( \overbar{\mb{B}}_r^\intercal \dot{\mb{B}}_r + \overbar{\mb{B}}_i^\intercal \dot{\mb{B}}_i \right) = \mathrm{Tr}\left( \overbar{\mb{A}}_r^\intercal \dot{\mb{A}}_r + \overbar{\mb{A}}_i^\intercal \dot{\mb{A}}_i \right),
\end{equation}
where we used the trace identity that $\tr{\mb{A}} = \tr{\mb{A}^\intercal}$.
\begin{sidewaystable}
% \begin{table}
\centering
\caption{RAD trace identity for different matrix functions.}
\begin{tabular}{ll}
\toprule
& \multicolumn{1}{l}{Formula} \\
\midrule
\multicolumn{1}{l}{Real differentiable} & $\tr{\overbar{\mb{B}}^\intercal \dot{\mb{B}}} = \tr{\overbar{\mb{A}}^\intercal \dot{\mb{A}}}$ (\citet{Giles2008})\\
\midrule
\multicolumn{1}{l}{imaginary analytic} & $\tr{\overbar{\mb{B}}^* \dot{\mb{B}}} = \tr{\overbar{\mb{A}}^* \dot{\mb{A}}}$ (\citet{He2023}) \\
\midrule
\multicolumn{1}{l}{imaginary multivariate differentiable} & \shortstack[l]{
$\tr{\overbar{\mb{B}}^* \dot{\mb{B}} + {(\overbar{\mb{B}}^* \dot{\mb{B}})}^*} = \tr{\overbar{\mb{A}}^* \dot{\mb{A}} + {(\overbar{\mb{A}}^* \dot{\mb{A}})}^*}$ (\citet{roberts2020a}) \\[1ex]
$\tr{ \overbar{\mb{B}}_r^\intercal \dot{\mb{B}}_r + \overbar{\mb{B}}_i^\intercal \dot{\mb{B}}_i} = \tr{\overbar{\mb{A}}_r^\intercal \dot{\mb{A}}_r + \overbar{\mb{A}}_i^\intercal \dot{\mb{A}}_i}$ (current paper)
} \\
\bottomrule
\end{tabular}\label{tab:dot_product_id}
% \end{table}
\end{sidewaystable}

\section{RAD formulae in the literature}
Giles~\cite{Giles2008} proposed the following RAD formula for singular value derivative 
\begin{equation}
\overbar{\mb{A}} = \mb{U}\overbar{\mb{S}}\mb{V}^\intercal.
\label{eq:giles_formula}
\end{equation}
This formula works with real-valued inputs and is only for the singular value derivative.
For the singular vector, an approach was shown in~\cite{Giles2008a} but no concrete implementable formula was proposed.
Townsend~\cite{townsend2016a} proposed the following RAD formula for singular variables derivative 
\begin{equation}
\begin{aligned}
&\overbar{\mb{A}} = [\mb{U}(\mb{F} \circ [\mb{U}^\intercal \overbar{\mb{U}} - \overbar{\mb{U}}^\intercal \mb{U}])\mb{S} + (\mb{I}_m - \mb{U}\mb{U}^\intercal)\overbar{\mb{U}}\mb{S}^{-1}]\mb{V}^\intercal \\
&+ \mb{U}(\mb{I}_k \circ \overbar{\mb{S}})\mb{V}^\intercal +\mb{U}[\mb{S}(\mb{F} \circ [\mb{V}^\intercal \overbar{\mb{V}} - \overbar{\mb{V}}^\intercal \mb{V}])\mb{V}^\intercal + \mb{S}^{-1}\overbar{\mb{V}}^\intercal(\mb{I}_n - \mb{V}\mb{V}^\intercal)],\\
\end{aligned}
\label{eq:townsend_formula}
\end{equation}
where $\mb{A}$ is an $m \times n$ matrix of rank $k \le \min({m,n})$, $\mb{U}$ is $m \times k$, $\mb{S}$ is $k \times k$, $\mb{V}$ is $n \times k$, $\mb{U}^\intercal \mb{U} = \mb{V}^\intercal \mb{V} = \mb{I}_k$ (the identity matrix of dimension $k$), and $\mb{F}$ is defined as 
\begin{equation}
\mathbf{F}_{ij} =
\begin{cases} 
\frac{1}{s_j^2 - s_i^2} & \text{if } i \neq j, \\ 
0 & \text{if } i = j,
\end{cases}
\label{eq:F}
\end{equation}
where $s$ is each entry of the singular values matrix $\mb{S}$.
Also, $\circ$ denotes the Haddamard product.
This formula works with truncated SVD and also with only real-valued inputs and outputs.
Wan and Zhang~\cite{wan2019a} proposed the following RAD formula 
\begin{equation}
\begin{aligned}
\overbar{\mb{A}} = \f{1}{2}\biggl(2\mb{U} \overbar{\mb{S}} \mb{V}^* + \mb{U}(\mb{J}+\mb{J}^*)\mb{S}\mb{V}^* + \mb{U}\mb{S}(\mb{K}+\mb{K}^* )\mb{V}^* + \f{1}{2} \mb{U} \mb{S}^{-1}(\mb{L}^*-\mb{L})\mb{V}^* \\
+ 2(\mb{I}-\mb{U}\mb{U}^*)\overbar{\mb{U}}\mb{S}^{-1}\mb{V}^* + 2\mb{U}\mb{S}^{-1} \overbar{\mb{V}}^*(1 - \mb{V} \mb{V}^*)\biggl),
\end{aligned}
\label{eq:wan_and_zhang}
\end{equation}
where $\mb{J} = \mb{F} \circ (2\mb{U}^* \overbar{\mb{U}})$, $\mb{K} = \mb{F} \circ (2\mb{V}^* \overbar{\mb{V}})$, $\mb{L} = \mb{I} \circ (2\mb{V}^* \overbar{\mb{V}})$, $\mb{F}$ is same as in \Cref{eq:F}, $\mb{I}$ is the identity matrix.
The formula proposed works with imaginary valued inputs and gives one single derivative value as opposed to the imaginary components shown in te current manuscript (such as $\d f_r / \d \mb{A}_r)$.
The details on this Wiritinger derivative are shown in~\Cref{sec:wiritinger}.

Seeger et al~\cite{seeger2017a} proposed the following RAD formula for singular variables derivative 
\begin{equation}
\begin{aligned}
\overline{\mb{A}} &= \mb{U}^\intercal \big( \mb{G}_2 \mb{V} + \Lambda^{-1} \overline{\mb{V}} \big),\\
\mb{G}_2 &= \mb{G}_1 \circ \mb{E} \Lambda - \big( \Lambda^{-1} \overline{\mb{V}} \mb{V}^\intercal \circ \mb{I} \big),\\
\mb{G}_1 &= \overline{\mb{U}}\mb{U}^\intercal + \Lambda^{-1} \overline{\mb{V}} \mb{V}^\intercal \mb{A},\\
\mb{E}_{i,j} &=
\begin{cases}
\frac{1}{\lambda_j - \lambda_i}, & i \neq j, \\
0, & i = j,
\end{cases}
\quad
h(t) = \max(|t|, \epsilon) \, \text{sgn}(t),
\end{aligned}
\label{eq:seeger}
\end{equation}
where the singular value decomposition (SVD) in Seeger et al~\cite{seeger2017a} is defined as
\begin{equation}
\mb{A} = \mb{U}^\intercal \Lambda \mb{V}, \quad \Lambda = \text{diag}(\lambda) \in \mathbb{R}^{m \times m}, \quad \mb{U} \in \mathbb{R}^{m \times m}, \quad \mb{V} \in \mathbb{R}^{m \times n},
\end{equation}
and further details on the computation of thin SVD can be found in Seeger et al~\cite{seeger2017a}.
The formula in \Cref{eq:seeger} works with real valued inputs.

\section{Wiritinger derivative}\label{sec:wiritinger}
In the current manuscript just as in He et al.~\cite{He2023}, the derivatives in the imaginary domain are expressed in terms of their real and imaginary components as $\d f_r / \d \mb{A}_r$ and such.
It may therefore be useful to the readers to see one single derivative $\d f/ \d \mb{A}$ for theimaginary valued case.
This can be achieved by using the Wiritinger derivative.
Information on this derivative can be found in the book by Hjorunges~\cite{hjorungnes2011a}.
The derivative is  
\begin{equation}
\f{\d f}{\d \mb{A}} = \f{1}{2} \bigl( \f{\d f_r}{\d \mb{A}_r} + i\f{\d f_i}{\d \mb{A}_r} - i\f{\d f_r}{\d \mb{A}_i} + \f{\d f_i}{\d \mb{A}_i} \bigl),
\end{equation}
where $i$ is iota, the imaginary number $\sqrt{-1}$.

For instance, take the singular value derivative as shown in \Cref{eq:Ar_bar_sigma_appendix,eq:Ai_bar_sigma_appendix}.
These go into the Wiritinger derivative as shown below
\begin{equation}
\f{\d \sigma}{\d \mb{A}} = \f{1}{2} \biggl[\f{\partial \sigma_r}{\partial \mb{A}_r} - i\f{\partial \sigma_r}{\partial \mb{A}_i}\biggl].
\label{eq:wiritinger_singular_var_appendix}
\end{equation}
The imaginary terms for the singular value have been removed in \Cref{eq:wiritinger_singular_var_appendix} because the singular value is purely real.
Before we proceed, consider the following equation:
\begin{equation}
\mb{u}\mb{v}^* = (\mb{u}_r + i \mb{u}_i) (\mb{v}_r^\intercal - i\mb{v}_i^\intercal)
= (\mb{u}_r \mb{v}_r^\intercal + \mb{u}_i \mb{v}_i^\intercal) + i(\mb{u}_i \mb{v}_r^\intercal - \mb{u}_r \mb{v}_i^\intercal).
\label{eq:uv*_appendix}
\end{equation}
We will then use \Cref{eq:Ar_bar_sigma_appendix} and \Cref{eq:Ai_bar_sigma_appendix} to expand \Cref{eq:wiritinger_singular_var_appendix}. 
After doing so, we obtain:
\begin{equation}
\f{\d \sigma}{\d \mb{A}} = \f{1}{2} \biggl[(\mb{u}_r \mb{v}_r^\intercal + \mb{u}_i \mb{v}_i^\intercal) - i(\mb{u}_i \mb{v}_r^\intercal - \mb{u}_r \mb{v}_i^\intercal)\biggl].
\label{eq:dsigmadA_form_2_appendix}
\end{equation}
Upon careful observation, one can notice that the term square in brackets of Eq.~\eqref{eq:dsigmadA_form_2_appendix} is simply the imaginary conjugate of Eq.~\eqref{eq:uv*_appendix}. 
Thus, we have
\begin{equation}
\f{\d \sigma}{\d \mb{A}} = \f{1}{2} {(\mb{u} \mb{v}^*)}^c,
\end{equation}
where $c$ denotes imaginary conjugate operation. 
In this way, we can compute the derivative of any one particular singular value with respect to its imaginary matrix $\mb{A}$.
The formula however must not be confused with the formula in the real case.
This is because the Wiritinger derivative is only applied in the imaginary plane treating the imaginary variable and its imaginary conjugate variable as two independent variables~\cite{hjorungnes2007a}.
Similarly, we can compute the total derivative from the imaginary parts in \Cref{eq:dfdA_components_SVD_Adjoint,eq:get_dfdA_from B,eq:get_dfdA_from_C}.

\section{Truss and finite element equations}\label{sec:truss}
\hlo{The governing equations for static response design sensitivity analysis is}
\begin{equation}
    \begin{aligned}
        \mb{K}\mb{d} = \mb{f}, 
    \end{aligned}
    \label{eq:struct_gov_eq}
\end{equation} % [x] TODO HS-: should this be reduced? (M.6)?
%RK-HS-: Yes Dr. He. 
\hlo{where $\mb{K} \in \mathbb{R}^{n_b \times n_b}$ is the global stiffness matrix, $\mb{f} \in \mathbb{R}^{n_b}$ is the vector of applied forces, $\mb{d}$ is the vector of displacements in the static condition and $n_b$ is the number of elements/bars in the truss.} % [x] TODO HS-: switch to different var for displacement? Maybe w, u? Make sure the dimension n is not overloaded...
%RK-HS-: w and u have been used. I replaced n with n_b as before for consistency.

\hlo{Next, we must discuss the finite element equations for}~\Cref{eq:struct_gov_eq}.
\hlo{Consider a 3D truss element with nodes $i$ and $j$, having coordinates $(x_i, y_i, z_i)$ and $(x_j, y_j, z_j)$ respectively. The element has cross-sectional area $b$, Young's modulus $E$ and length $L$.}

\hlo{The direction cosines of the element are}
\begin{equation}
l = \frac{x_j - x_i}{L}, \quad m = \frac{y_j - y_i}{L}, \quad n = \frac{z_j - z_i}{L}.
\end{equation}
\hlo{The local stiffness matrix in global coordinates (6$\times$6) is}
\begin{equation}
\mathbf{K}^{(e)}_{\text{local}} = \frac{bE}{L}
\begin{bmatrix}
l^2 & lm & ln & -l^2 & -lm & -ln \\
lm & m^2 & mn & -lm & -m^2 & -mn \\
ln & mn & n^2 & -ln & -mn & -n^2 \\
-l^2 & -lm & -ln & l^2 & lm & ln \\
-lm & -m^2 & -mn & lm & m^2 & mn \\
-ln & -mn & -n^2 & ln & mn & n^2
\end{bmatrix},
\label{eq:local_stiff} % [x] TODO HS-: capitalize k^e? It is a matrix.
\end{equation}
\hlo{where $bE/L$ is basically the stiffness of the bar.
The global stiffness matrix $\mathbf{K}$ is assembled by adding contributions of all element stiffness matrices $\mathbf{K}^{(e)}_{\text{local}}$ to their appropriate global degrees of freedom (DOFs).}

\hlo{For each element, we map the local DOFs to global ones. Suppose an element connects nodes $i$ and $j$, the global DOF vector is:}
\begin{equation}
\mb{u}^{(e)} = [u_{i_x}, u_{i_y}, u_{i_z}, u_{j_x}, u_{j_y}, u_{j_z}]^T.
\end{equation}
\hlo{Next, we assemble each element's stiffness matrix into the global matrix using these DOF indices}
\begin{equation}
\mathbf{K} = \sum_{e=1}^{n_{\text{elem}}} \mathbf{A}_e^\intercal \mathbf{K}^{(e)}_{\text{local}} \mathbf{A}_e,
\label{eq:local_global}
\end{equation}
\hlo{where $\mathbf{A}_e$ is a Boolean assembly matrix that maps local to global DOFs. % [x] TODO HS-: what is a Boolean assembly matrix? Either define it here or add ref.
In the finite element method (FEM), a ``Boolean assembly matrix'' is used to map local DOFs from individual elements to the global system. It is called ``Boolean'' because its entries are either 0 or 1, indicating whether a particular local DOF contributes to a global DOF.

Consider a system with $a$ global DOF and an element with $b$ local DOFs.
When assembling the global stiffness matrix $\mb{K}$ from local element matrices $\mb{K}^{(e)}_{\text{local}}$, the Boolean matrix is used as in}~\Cref{eq:local_global} \hlo{where each row of $\mb{A}_e$ has a single entry of 1 corresponding to the global index of the respective local DOF. This facilitates correct placement of local element contributions into the global system matrix or vector.}

% [x] TODO HS-: mathbf{K}_{\text{global}} or just K?
\hlo{Application of boundary conditions involves modifying $\mathbf{K}$ and $\mathbf{f}$ to enforce prescribed displacements, zeroing rows and columns of fixed DOFs and inserting 1 along the diagonal, and setting corresponding entries in $\mathbf{f}$ to the prescribed values.
The reduction for $\mb{f}$ is done by identifying an index set of free DOFs and extracting the subvector $ \mathbf{f}_{\text{red}}$, which is the reduced vector of $\mb{f}$.}

% [x] TODO HS-: \mathbf{f}_{\text{mod}}, or just f?
% [x] TODO HS-: talk about how to assemble f before the global elastic equation?
% [x] TODO HS-: avoid using {} for vector...

\hlo{Finally, we solve the modified system}
\begin{equation}
\mathbf{K}_{\text{red}} \mathbf{d}_{\text{red}} = \mathbf{f}_{\text{red}},
\end{equation}
\hlo{where $\mb{K}_{\text{red}}$ is the reduced stiffness matrix, $\mathbf{d}_{\text{red}}$ is the vector of displacements that with only the nodes we solve for and $\mathbf{f}_{\text{red}}$ is the reduced vector of forces applied to only the nodes under consideration.
We refer the readers to Chandrupatla et al.[~\cite{Chandrupatla2021},Chapter 4] which throws light on the three-dimensional finite element formulation for trusses.

It is noted at this stage that we perform SVD on the inverse this reduced stiffness matrix as}
\begin{equation}
        \mb{G} = \mb{K}_{\text{red}}^{-1},
    \label{eq:transfer_func_mat}
\end{equation}
\hlo{where $\mb{G}$ is the transformation matrix.}

\hlo{Once $\mathbf{d}_{\text{mod}}$ is obtained, element forces can be calculated using:}
\begin{equation}
f^{(e)} = \frac{bE}{L} \begin{bmatrix} -l & -m & -n & l & m & n \end{bmatrix} \mb{d}^{(e)}_{\text{red}},
\end{equation}
\hlo{where $f$ is the internal force for the element $e$.}

\hlo{For the response of the system, SVD of the matrix} $\mb{G}$ \hlo{is sought as shown in}~\Cref{eq:transfer_func_mat} \hlo{which gives the singular variables subject to design conditions in $\mb{b}$.}
\hlo{In the study of Ersoy and Mugan~\cite{ersoy2002a}, the design parameter $b_1$ was found to be the most effective parameter in order to reduce displacements, or in other words, increase the stiffness of the structure.}
\hlo{In the following part of this section, we formulate RAD-based equations which includes the RAD form for the singular value sensitivity from}~\Cref{eq:RAD_formula_singular_sens_real}.

\hlo{The path of perturbations is first studied.
The design variables affect the stiffness matrix $\mb{K}_{\text{red}}$.
$\mb{K}_{\text{red}}$ in turn affects $\mb{G}$, the transfer function matrix.
The SVD of $\mb{G}$ gives the first singular value $\sigma_1$, whose derivative was found to be crucial in minimization of the total static response magnitude~\cite{ersoy2002a}.}

\hlo{The derivative of the first singular value with respect to the matrix $\mb{G}$ is} % [x] TODO HS-: avoid "given as" check. Also directly refer wheere you derive the equations.
\begin{equation}
    \f{\d \sigma_1}{\d \mb{G}} = \mb{u}_1 \mb{v}_1^\intercal,
    \label{eq:sing_val_der_G}
\end{equation}
\hlo{where $\mb{u},\mb{v}$ is the first pair of singular vectors corresponding to $\sigma_1$.
This was derived earlier in}~\Cref{sec:Singular_Val_Sens_RAD_Derivation}\hlo{ in}~\Cref{eq:appendix_RAD_sigma}.
\hlo{The derivative is given in}~\Cref{eq:RAD_formula_singular_sens_real}.
\hlo{Next, the RAD form for}~\Cref{eq:transfer_func_mat} \hlo{is given by}
\begin{equation}
    \overbar{\mb{K}}_{\text{red}} = \mb{K}_{\text{red}}^{-\intercal} \overbar{\mb{G}} \mb{K}_{\text{red}}^{-\intercal},
\end{equation}
\hlo{in which upon seeding $\sigma_1$ as a reverse seed, we get}
\begin{equation}
    \f{\d \sigma_1}{\d \mb{K}_{\text{red}}} = \mb{K}^{-\intercal}_{\text{red}} \mb{u}_1 \mb{v}_1^\intercal \mb{K}^{-\intercal}_{\text{red}},
\end{equation}
\hlo{which is the derivative we require before we flatten the derivative matrix in LHS into a vector.}
\hlo{For the derivation of the RAD form for inverse of a real-valued matrix, we refer the readers to Giles}~\cite{Giles2008a}.

\hlo{The final step is to obtain $\d \mb{K}_{\text{red}} / \d \mb{b}$.}
\hlo{It is evident from}~\Cref{eq:local_stiff} \hlo{that the derivative of $\d \mb{K}_{\text{red}} / \d \mb{b}$ is a tensor.
$\d \mb{K}_{\text{red}} / \d b_1$ for instance is simply}
\begin{equation}
    \f{\d \mathbf{K}^{(e)}_{\text{red}}}{\d b_1} = \frac{E}{L}
    \begin{bmatrix}
    l^2 & lm & ln & -l^2 & -lm & -ln \\
    lm & m^2 & mn & -lm & -m^2 & -mn \\
    ln & mn & n^2 & -ln & -mn & -n^2 \\
    -l^2 & -lm & -ln & l^2 & lm & ln \\
    -lm & -m^2 & -mn & lm & m^2 & mn \\
    -ln & -mn & -n^2 & ln & mn & n^2
    \end{bmatrix},
    \label{eq:dK_b}
\end{equation}
% [x] TODO HS-: check this equation...
\hlo{which is simply a matrix of all constants with $e=1$.
This matrix does not change for the other areas in $\mb{b}$.
Each matrix in this tensor is then flattened into a vector.
This tensor is then simply plugged into}~\Cref{eq:der_RAD_opt} \hlo{to complete the derivative formula for $\d \sigma_1 / \d \mb{b}$.}

\end{document}